\documentclass[10pt,twoside]{amsart}
\usepackage{amsfonts,amssymb,amsmath}

\def\YEAR{\year}\newcount\VOL\VOL=\YEAR\advance\VOL by-1995
\def\firstpage{1}\def\lastpage{1000}
\def\received{}\def\revised{}
\def\communicated{}

\makeatletter
\def\magnification{\afterassignment\m@g\count@}
\def\m@g{\mag=\count@\hsize6.5truein\vsize8.9truein\dimen\footins8truein}
\makeatother

\font\eightrm=cmr8
\font\caps=cmcsc10                    
\font\Caps=cmcsc10 scaled \magstep1   

\makeatletter
\@specialpagefalse\headheight=8.5pt
\def\DocMath{}
\renewcommand{\@evenhead}{%
    \ifnum\thepage>\lastpage\rlap{\thepage}\hfill%
    \else\rlap{\thepage}\slshape\leftmark\hfill{\caps\SAuthor}\hfill\fi}%
\renewcommand{\@oddhead}{%
    \ifnum\thepage=\firstpage{\DocMath\hfill\llap{\thepage}}%
    \else{\slshape\rightmark}\hfill{\caps\STitle}\hfill\llap{\thepage}\fi}%
\makeatother

\def\TSkip{\bigskip}
\newbox\TheTitle{\obeylines\gdef\GetTitle #1
\ShortTitle  #2
\SubTitle    #3
\Author      #4
\ShortAuthor #5
\EndTitle
{\setbox\TheTitle=\vbox{\baselineskip=20pt\let\par=\cr\obeylines%
\halign{\centerline{\Caps##}\cr\noalign{\medskip}\cr#1\cr}}%
        \copy\TheTitle\TSkip\TSkip%
\def\next{#2}\ifx\next\empty\gdef\STitle{#1}\else\gdef\STitle{#2}\fi%
\def\next{#3}\ifx\next\empty%
    \else\setbox\TheTitle=\vbox{\baselineskip=20pt\let\par=\cr\obeylines%
    \halign{\centerline{\caps##} #3\cr}}\copy\TheTitle\TSkip\TSkip\fi%
\centerline{\caps #4}\TSkip\TSkip%
\def\next{#5}\ifx\next\empty\gdef\SAuthor{#4}\else\gdef\SAuthor{#5}\fi%
\ifx\received\empty\relax
    \else\centerline{\eightrm Received: \received}\fi%
\ifx\revised\empty\TSkip%
    \else\centerline{\eightrm Revised: \revised}\TSkip\fi%
\ifx\communicated\empty\relax
    \else\centerline{\eightrm Communicated by \communicated}\fi\TSkip\TSkip%
\catcode'015=5}}\def\Title{\obeylines\GetTitle}
\def\Abstract{\begingroup\narrower
    \parskip=\medskipamount\parindent=0pt{\caps Abstract. }}
\def\EndAbstract{\par\endgroup\TSkip}

\long\def\MSC#1\EndMSC{\def\arg{#1}\ifx\arg\empty\relax\else
     {\par\narrower\noindent%
     2000 Mathematics Subject Classification: #1\par}\fi}

\long\def\KEY#1\EndKEY{\def\arg{#1}\ifx\arg\empty\relax\else
        {\par\narrower\noindent Keywords and Phrases: #1\par}\fi\TSkip}

\newbox\TheAdd\def\Addresses{\vfill\copy\TheAdd\vfill
    \ifodd\number\lastpage\vfill\eject\phantom{.}\vfill\eject\fi}
{\obeylines\gdef\GetAddress #1
\Address #2 
\Address #3
\Address #4
\EndAddress
{\def\xs{4.3truecm}\parindent=0pt
\setbox0=\vtop{{\obeylines\hsize=\xs#1\par}}\def\next{#2}
\ifx\next\empty 
     \setbox\TheAdd=\hbox to\hsize{\hfill\copy0\hfill}
\else\setbox1=\vtop{{\obeylines\hsize=\xs#2\par}}\def\next{#3}
\ifx\next\empty 
     \setbox\TheAdd=\hbox to\hsize{\hfill\copy0\hfill\copy1\hfill}
\else\setbox2=\vtop{{\obeylines\hsize=\xs#3\par}}\def\next{#4}
\ifx\next\empty\ 
     \setbox\TheAdd=\vtop{\hbox to\hsize{\hfill\copy0\hfill\copy1\hfill}
                \vskip20pt\hbox to\hsize{\hfill\copy2\hfill}}
\else\setbox3=\vtop{{\obeylines\hsize=\xs#4\par}}
     \setbox\TheAdd=\vtop{\hbox to\hsize{\hfill\copy0\hfill\copy1\hfill}
                \vskip20pt\hbox to\hsize{\hfill\copy2\hfill\copy3\hfill}}
\fi\fi\fi\catcode'015=5}}\gdef\Address{\obeylines\GetAddress}

\begin{document}
\Title Non-Hausdorff groupoids, proper actions and $K$-theory
\ShortTitle Non-Hausdorff groupoids
\SubTitle   
\Author Jean-Louis Tu
\ShortAuthor 
\EndTitle
\Abstract Let $G$ be a (not necessarily Hausdorff) locally compact groupoid.
We introduce a notion of properness for $G$,
which is invariant under Morita-equivalence.
We show that any generalized morphism between two locally
compact groupoids which satisfies some properness conditions
induces a $C^*$-correspondence from
$C^*_r(G_2)$ to $C^*_r(G_1)$, and thus two Morita equivalent
groupoids have Morita-equivalent $C^*$-algebras.
\EndAbstract
\MSC 22A22 (Primary); 46L05, 46L80, 54D35 (Secondary).
\EndMSC
\KEY groupoid, $C^*$-algebra, $K$-theory.
\EndKEY
\Address
Jean-Louis Tu
University Paris VI
Institut de Math\'ematiques
175, rue du Chevaleret
75013 Paris, France.
{\tt tu@math.jussieu.fr}
\Address
\Address
\Address
\EndAddress

\newtheorem{theo}{Theorem}[section]
\newtheorem{prop}[theo]{Proposition}
\newtheorem{lem}[theo]{Lemma}
\newtheorem{coro}[theo]{Corollary}

\newtheorem{defi}[theo]{Definition}
\newtheorem{example}[theo]{Example}
\newtheorem{examples}[theo]{Examples}
\newtheorem{rem}[theo]{Remark}
\newtheorem{rems}[theo]{Remarks}


\newcommand\Z{{\mathbb{Z}}}
\newcommand\N{{\mathbb{N}}}
\newcommand\R{{\mathbb{R}}}
\newcommand\C{{\mathbb{C}}}

\def\rond{\mathaccent"7017}


\newcommand\pf{\begin{proof}}

\newcommand\pfend{\end{proof}}

\section*{Introduction}
Very often, groupoids that appear in geometry, such as
holonomy groupoids of foliations, groupoids of inverse semigroups
\cite{pat,ks02} and the indicial algebra of a manifold with
corners \cite{lm01} are not Hausdorff.
It is thus necessary to extend various basic notions to this
broader setting, such as proper action and Morita equivalence.
We also show that a generalized morphism from $G_2$ to $G_1$
satisfying certain properness conditions induces an element
of $KK(C^*_r(G_2),C^*_r(G_1))$.
%
%
%
\par\bigskip


In Section~\ref{sec:proper groupoids}, we introduce the notion
of proper groupoids and show that it is invariant under Morita-equivalence.

Section~\ref{sec:topological construction} is a technical part of the
paper in which from every locally compact topological space $X$ is
canonically constructed a locally compact Hausdorff space ${\mathcal{H}}X$
in which
$X$ is (not continuously) embedded.
When $G$ is a groupoid (locally compact, with Haar system, such
that $G^{(0)}$ is Hausdorff), the closure $X'$ of $G^{(0)}$ in
${\mathcal{H}}G$ is endowed with a continuous action of $G$
and plays an important technical r\^ole.

In Section~\ref{sec:haar systems}  we review basic properties
of locally compact groupoids with Haar system and technical
tools that are used later.

In Section~\ref{sec:Hilbert module of proper groupoid} we construct,
using tools of Section~\ref{sec:topological construction},
a canonical $C^*_r(G)$-Hilbert module ${\mathcal{E}}(G)$
for every (locally compact...) proper groupoid $G$. If $G^{(0)}/G$ is
compact, then there exists a projection $p\in C^*_r(G)$ such that
${\mathcal{E}}(G)$ is isomorphic to $pC^*_r(G)$. The
projection $p$ is given by $p(g)=(c(s(g))c(r(g)))^{1/2}$,
where $c\colon G^{(0)}\to \R_+$ is a ``cutoff'' function
(Section~\ref{sec:cutoff}).
Contrary to the Hausdorff case, the function $c$ is not continuous,
but it is the restriction to $G^{(0)}$ of a continuous map
$X'\to \R_+$ (see above for the definition of $X'$).

In Section~\ref{sec:correspondences}, we examine the
question of naturality $G\mapsto C^*_r(G)$. Recall
that if $f\colon X\to Y$ is a continuous map between two locally
compact spaces, then $f$ induces a map from $C_0(Y)$ to $C_0(X)$
if and only if $f$ is proper. When $G_1$ and $G_2$ are groups,
a morphism $f\colon G_1\to G_2$ does not induce a map
$C^*_r(G_2)\to C^*_r(G_1)$ (when $G_1\subset G_2$ is an inclusion
of discrete groups there is a map in the other direction).
When $f\colon G_1\to G_2$ is a groupoid morphism, we cannot expect
to get more than a $C^*$-correspondence from $C^*_r(G_2)$ to
$C^*_r(G_1)$ when $f$ satisfies certain properness assumptions:
this was done in the Hausdorff situation by Macho-Stadler and
O'Uchi (\cite[Theorem 2.1]{mo}, see also \cite{lan,mrc,sta}),
but the formulation of their
theorem is somewhat complicated. In this paper, as a corollary of
Theorem~\ref{thm:correspondence}, we get that (in the Hausdorff
situation), if the restriction of $f$ to $(G_1)_K^K$ is proper
for each compact set $K\subset (G_1)^{(0)}$ then $f$ induces
a correspondence ${\mathcal{E}}_f$ from $C^*_r(G_2)$ to
$C^*_r(G_1)$. In fact we construct a $C^*$-correspondence
out of any groupoid generalized morphism (\cite{hs,leg})
which satisfies some properness conditions. As a corollary,
if $G_1$ and $G_2$ are Morita equivalent then $C^*_r(G_1)$
and $C^*_r(G_2)$ are Morita-equivalent $C^*$-algebras.
\par\medskip

Finally, let us add that our original motivation was to
extend Baum, Connes and Higson's construction of the assembly
map $\mu$ to non-Hausdorff groupoids; however, we couldn't
prove $\mu$ to be an isomorphism in any non-trivial case.

%
%

\section{Preliminaries}
\subsection{Groupoids}
Throughout, we will assume that the reader is familiar with basic
definitions about groupoids
(see~\cite{ren,pat}). If $G$ is a groupoid, we denote by $G^{(0)}$ its
set of units and by $r\colon G\to G^{(0)}$ and
$s\colon G\to G^{(0)}$ its range and
source maps respectively. We will use notations such as
$G_x=s^{-1}(x)$, $G^y=r^{-1}(y)$, $G_x^y=G_x\cap G^y$.
Recall that a topological groupoid is said
to be {\em \'{e}tale}\/ if $r$ (and $s$) are local homeomorphisms.
\par\medskip

For all sets $X$, $Y$, $T$ and all maps $f\colon X\to T$
and $g\colon Y\to T$, we denote by $X\times_{f,g}Y$, or by
$X\times_T Y$ if there is no ambiguity, the
set $\{(x,y)\in X\times Y\vert\; f(x)=g(y)\}$.\par\medskip

Recall that a (right) action of $G$ on a set $Z$ is given by
\begin{itemize}
\item[{\rm (a)}] a (``momentum'') map $p\colon Z\to G^{(0)}$;
\item[{\rm (b)}] a map $Z\times_{p,r}G\to Z$, denoted
by $(z,g)\mapsto zg$
\end{itemize}
with the following properties:
\begin{itemize}
\item[{\rm (i)}] $p(zg)=s(g)$ for all $(z,g)\in Z\times_{p,r}G$;
\item[{\rm (ii)}] $z(gh)=(zg)h$ whenever $p(z)=r(g)$ and $s(g)=r(h)$;
\item[{\rm (iii)}] $zp(z)=z$ for all $z\in Z$.
\end{itemize}
Then the
crossed-product $Z\rtimes G$ is the subgroupoid of $(Z\times Z)
\times G$ consisting of elements $(z,z',g)$ such that $z'=zg$.
Since the map $Z\rtimes G\to Z\times G$ given by
$(z,z',g)\mapsto (z,g)$ is
injective, the groupoid $Z\rtimes G$ can also be considered as
a subspace of $Z\times G$, and this is what we will do most of the time.

\subsection{Locally compact spaces}
A topological space $X$ is said to be quasi-compact if every open
cover of $X$ admits a finite sub-cover. A space is compact if it is
quasi-compact and Hausdorff.
Let us recall a few basic facts about locally compact spaces.
\begin{defi}
A topological space $X$ is said to be locally compact if every point
$x\in X$ has a compact neighborhood.
\end{defi}

In particular, $X$ is locally Hausdorff,
thus every singleton subset of $X$ is
closed. Moreover, the diagonal in $X\times X$ is locally closed.

\begin{prop}
Let $X$ be a locally compact space. Then every locally closed subspace of
$X$ is locally compact.
\end{prop}

Recall that $A\subset X$ is locally closed if for every $a\in A$,
there exists a neighborhood $V$ of $a$ in $X$ such that $V\cap A$ is
closed in $V$. Then $A$ is locally closed if and only if it is of the
form $U\cap F$, with $U$ open and $F$ closed.

\begin{prop}
Let $X$ be a locally compact space. The following are equivalent:
\begin{itemize}
\item[\rm{(i)}] there exists a sequence $(K_n)$ of compact subspaces
  such that $X=\cup_{n\in \N} K_n$; 
\item[\rm{(ii)}] there exists a sequence $(K_n)$ of quasi-compact subspaces
  such that $X=\cup_{n\in \N} K_n$; 
\item[\rm{(iii)}] there exists a sequence $(K_n)$ of quasi-compact subspaces
  such that $X=\cup_{n\in \N} K_n$ and $K_n\subset \rond{K}_{n+1}$ for
  all $n\in \N$.
\end{itemize}
\end{prop}

Such a space will be called $\sigma$-compact.
\pf
(i)$\implies$(ii) is obvious. The implications
(ii)$\implies$(iii)$\implies$(i)
follow easily from the fact that for every quasi-compact subspace $K$,
there exists a finite family $(K_i)_{i\in I}$ of compact sets such
that $K\subset \cup_{i\in I} \rond{K_i}$. 
\pfend

\subsection{Proper maps}
\begin{prop}\cite[Th\'{e}or\`{e}me I.10.2.1]{bou}\label{prop:proper maps}
Let $X$ and $Y$ be two topological spaces, and $f\colon X\to Y$ a
continuous map. The following are equivalent:

\begin{itemize}
\item[\rm{(i)}] For every topological space $Z$,
  $f\times{\mathrm{Id}}_Z\colon X\times Z\to Y\times Z$ is closed;
\item[\rm{(ii)}] $f$ is closed and for every $y\in Y$, $f^{-1}(y)$ is
  quasi-compact.
\end{itemize}
\end{prop}

A map which satisfies the equivalent properties of
Proposition~\ref{prop:proper maps} is said to be {\it proper}.

\begin{prop}\cite[Proposition I.10.2.6]{bou}
\label{prop:proper implies inverse image quasi-compact}
Let $X$ and $Y$ be two topological spaces and let $f\colon X\to Y$ be
a proper map. Then for every quasi-compact subspace $K$ of $Y$,
$f^{-1}(K)$ is quasi-compact.
\end{prop}

\begin{prop}\label{prop:characterization of proper maps}
Let $X$ and $Y$ be two topological spaces and let $f\colon X\to Y$ be
a continuous map. Suppose $Y$ is locally compact, then the following are
equivalent:
\begin{itemize}
\item[{\rm{(i)}}] $f$ is proper;
\item[{\rm{(ii)}}] for every quasi-compact subspace $K$ of $Y$,
  $f^{-1}(K)$ is quasi-compact;
\item[{\rm{(iii)}}] for every compact subspace $K$ of $Y$,
  $f^{-1}(K)$ is quasi-compact;
\item[{\rm{(iv)}}] for every $y\in Y$, there exists a compact
  neighborhood $K_y$ of $y$ such that $f^{-1}(K_y)$ is quasi-compact.
\end{itemize}
\end{prop}

\pf
(i)$\implies$(ii) follows from Proposition~\ref{prop:proper implies
  inverse image quasi-compact}.
(ii)$\implies$(iii)$\implies$(iv) are
obvious. Let us show (iv)$\implies$(i).

Since $f^{-1}(y)$ is closed, it is clear that $f^{-1}(y)$ is
quasi-compact for all $y\in Y$. It
remains to prove that for every closed subspace $F\subset X$, $f(F)$ is
closed. Let $y\in \overline{f(F)}$.
Let $A=f^{-1}(K_y)$.
Then $A\cap F$ is quasi-compact, so $f(A\cap F)$ is
quasi-compact. As $f(A\cap F)\subset K_y$, it is closed in
$K_y$, i.e. $K_y\cap\overline{f(A\cap F)}=K_y\cap f(A\cap F)$.
We thus have
$y\in K_y\cap\overline{f(A\cap F)}=K_y\cap f(A\cap
F)\subset f(F)$. It follows that $f(F)$
is closed.
\pfend

\section{Proper groupoids and proper actions}\label{sec:proper groupoids}
\subsection{Locally compact groupoids}
\begin{defi}
A topological groupoid $G$ is said to be locally compact
(resp. $\sigma$-compact) if it is locally compact
(resp. $\sigma$-compact) as a topological space.
\end{defi}

\begin{rem}
The definition of a locally compact groupoid in \cite{pat} corresponds
to our definition of a locally compact, $\sigma$-compact groupoid with
Haar system whose unit space is Hausdorff, thanks to
Propositions~\ref{prop:G0 locally compact} and~\ref{prop:Gx Hausdorff}.
\end{rem}

\begin{example}\label{ex:NH gpd}
Let $\Gamma$ be a discrete group, $H$ a closed
normal subgroup and let $G$ be the bundle of groups over $[0,1]$
such that $G_0=\Gamma$ and $G_t=\Gamma/H$ for all $t>0$.
We endow $G$ with the quotient topology of $\left([0,1]\times \Gamma
\right)/\left((0,1]\times H\right)$. Then $G$ is a non-Hausdorff
locally compact groupoid such that $(t,\bar\gamma)$ converges to
$(0,\gamma h)$ as $t\to 0$, for all $\gamma\in \Gamma$ and
$h\in H$.
\end{example}

\begin{example}\label{ex:NH foliation}
Let $\Gamma$ be a discrete group acting on a locally compact Hausdorff
space $X$, and let
$G=(X\times\Gamma)/\sim$, where $(x,\gamma)$ and $(x,\gamma')$ are
identified if their germs are equal, i.e. there exists a neighborhood
$V$ of $x$ such that $y\gamma=y\gamma'$ for all $y\in V$. Then
$G$ is locally compact, since the open sets $V_\gamma=
\{[(x,\gamma)]\vert\; x\in X\}$ are homeomorphic to $X$ and cover $G$.

Suppose that $X$ is a manifold, $M$ is a manifold such that
$\pi_1(M)=\Gamma$, $\tilde M$ is the universal cover of $M$ and
$V=(X\times \tilde M)/\Gamma$, then $V$ is foliated by
$\{[x,\tilde m]\vert\; \tilde m\in \tilde M\}$ and $G$ is the
restriction to a transversal of the holonomy groupoid of the above
foliation.
\end{example}

\begin{prop}\label{prop:G0 locally compact}
If $G$ is a locally compact groupoid, then $G^{(0)}$ is locally closed
in $G$, hence locally
compact. If furthermore $G$ is $\sigma$-compact, then $G^{(0)}$ is
$\sigma$-compact.
\end{prop}

\pf
Let $\Delta$ be the diagonal in $G\times G$. Since $G$ is locally
Hausdorff, $\Delta$ is locally closed. Then $G^{(0)}=(\mbox{Id},
r)^{-1}(\Delta)$ is locally closed in $G$.

Suppose that $G=\cup_{n\in \N}K_n$ with $K_n$ quasi-compact, then
$s(K_n)$ is quasi-compact and $G^{(0)}=\cup_{n\in \N}s(K_n)$.
\pfend

\begin{prop}
Let $Z$ a locally compact space and $G$ be a locally compact groupoid
acting on $Z$. Then the crossed-product $Z\rtimes G$ is locally compact.
\end{prop}

\pf
Let $p\colon Z\to G^{(0)}$ be the momentum map of the action of $G$.
From Proposition~\ref{prop:G0 locally compact}, the diagonal $\Delta\subset
G^{(0)}\times G^{(0)}$ is locally closed in $G^{(0)}\times G^{(0)}$,
hence $Z\rtimes G=(p,r)^{-1}(\Delta)$ is locally closed in $Z\times G$.
\pfend

Let $T$ be a space. Recall that there is a groupoid
$T\times T$ with unit space $T$, and product $(x,y)(y,z)=(x,z)$.

Let $G$ be a groupoid and $T$ be a space. Let $f\colon T\to G^{(0)}$,
and let $G[T]=\{(t',t,g)\in (T\times T)\times G\vert\; g\in
G_{f(t)}^{f(t')}\}$. Then $G[T]$ is a subgroupoid of $(T\times T)\times G$.

\begin{prop}
Let $G$ be a topological groupoid with $G^{(0)}$ locally Hausdorff,
$T$ a topological space and
$f\colon T\to G^{(0)}$ a continuous map. Then $G[T]$ is a locally
closed subgroupoid of $(T\times T)\times G$. In particular, if $T$ and
$G$ are locally compact, then $G[T]$ is locally compact.
\end{prop}

\pf
Let $F\subset T\times G^{(0)}$ be the graph of $f$. Then
$F=(f\times\mbox{Id})^{-1}(\Delta)$, where $\Delta$ is the diagonal in
$G^{(0)}\times G^{(0)}$, thus it is locally closed. Let
$\rho\colon (t',t,g)\mapsto (t',r(g))$ and
$\sigma\colon (t',t,g)\mapsto (t,s(g))$
be the range and source maps of $(T\times T)\times G$, then
$G[T]=(\rho,\sigma)^{-1}(F\times F)$ is locally closed.
\pfend

\begin{prop}\label{prop:Gx Hausdorff}
Let $G$ be a locally compact groupoid such that $G^{(0)}$ is
Hausdorff. Then for every $x\in G^{(0)}$, $G_x$ is Hausdorff.
\end{prop}

\pf
Let $Z=\{(g,h)\in G_x\times G_x\vert\; r(g)=r(h)\}$.
Let $\varphi\colon Z\to G$ defined by
$\varphi(g,h)=g^{-1}h$. Since $\{x\}$ is closed in $G$,
$\varphi^{-1}(x)$ is closed in $Z$, and since $G^{(0)}$ is Hausdorff,
$Z$ is closed in $G_x\times G_x$. It follows that $\varphi^{-1}(x)$,
which is the diagonal of $G_x\times G_x$, is closed in $G_x\times G_x$.
\pfend

\subsection{Proper groupoids}
\begin{defi}
A topological groupoid $G$ is said to be proper
if $(r,s)\colon G\to G^{(0)}\times G^{(0)}$ is proper.
\end{defi}

\begin{prop}\label{prop:proper groupoid}
Let $G$ be a topological groupoid such that $G^{(0)}$ is locally compact.
Consider the following assertions:
\begin{itemize}
\item[{\rm{(i)}}] $G$ is proper;
\item[{\rm{(ii)}}] $(r,s)$ is closed and for every $x\in G^{(0)}$,
  $G_x^x$ is quasi-compact;
\item[{\rm{(iii)}}] for all quasi-compact subspaces $K$ and $L$ of
  $G^{(0)}$, $G_K^L$ is quasi-compact;
\item[{\rm{(iii)'}}] for all compact subspaces $K$ and $L$ of
  $G^{(0)}$, $G_K^L$ is quasi-compact;
\item[{\rm{(iv)}}] for every quasi-compact subspace $K$ of $G^{(0)}$, $G_K^K$
  is quasi-compact;
\item[{\rm{(v)}}] $\forall x$, $y\in G^{(0)}$, $\exists K_x$, $L_y$
  compact neighborhoods of $x$ and $y$ such that $G_{K_x}^{L_y}$ is
  quasi-compact.
\end{itemize}
Then (i)$\iff$(ii)$\iff$(iii)$\iff$(iii)'$\iff$(v)$\implies$(iv).
If $G^{(0)}$ is Hausdorff, then (i)--(v) are equivalent.
\end{prop}

\pf
(i)$\iff$(ii) follows from Proposition~\ref{prop:proper maps}, and
from the fact that $G_x^x$ is homeomorphic to $G_x^y$ if
$G_x^y\ne\emptyset$.
(i)$\implies$(iii) and (v)$\implies$(i) follow
Proposition~\ref{prop:characterization of proper maps} and the formula
$G_K^L=(r,s)^{-1}(L\times K)$.
(iii)$\implies$(iii)'$\implies$(v) and (iii)$\implies$(iv) are obvious.
If $G^{(0)}$ is Hausdorff, then (iv)$\implies$(v) is obvious.
\pfend

Note that if $G=G^{(0)}$ is a non-Hausdorff topological space, then
$G$ is not proper (since $(r,s)$ is not closed), but satisfies
property (iv).

\begin{prop}\label{prop:pi open}
Let $G$ be a topological groupoid. If $r\colon G\to G^{(0)}$ is open
then the canonical mapping $\pi\colon G^{(0)}\to G^{(0)}/G$ is
open.
\end{prop}

\pf
Let $V\subset G^{(0)}$ be an open
subspace. If $r$ is open, then
$r(s^{-1}(V))=\pi^{-1}(\pi(V))$ is open. Therefore, $\pi(V)$ is open.
\pfend

\begin{prop}\label{prop:G0/G}
Let $G$ be a topological groupoid such that $G^{(0)}$ is locally
compact and $r\colon G\to G^{(0)}$ is open. Suppose that $(r,s)(G)$ is
locally closed in $G^{(0)}\times G^{(0)}$, then
$G^{(0)}/G$ is locally compact. Furthermore,
\begin{itemize}
\item[{\rm{(a)}}] if $G^{(0)}$
is $\sigma$-compact, then $G^{(0)}/G$ is $\sigma$-compact;
\item[{\rm{(b)}}] if $(r,s)(G)$ is closed (for instance if $G$ is
  proper), then $G^{(0)}/G$ is Hausdorff.
\end{itemize}
\end{prop}

\pf
Let $R=(r,s)(G)$. Let $\pi\colon G^{(0)}\to G^{(0)}/G$ be the canonical
mapping. By Proposition~\ref{prop:pi open}, $\pi$ is open, therefore
$G^{(0)}/G$ is locally quasi-compact. Let us show that it is locally
Hausdorff. Let $V$ be an open subspace of $G^{(0)}$ such that $(V\times
V)\cap R$ is closed in $V\times V$. Let $\Delta$ be the diagonal in
$\pi(V)\times \pi(V)$. Then $(\pi\times\pi)^{-1}(\Delta)=(V\times
V)\cap R$ is closed in $V\times V$. Since $\pi\times \pi\colon V\times
V\to \pi(V)\times\pi(V)$ is continuous open surjective, it follows
that $\Delta$ is closed in $\pi(V)\times\pi(V)$, hence $\pi(V)$ is
Hausdorff. This completes the proof that $G^{(0)}/G$ is locally
compact and of assertion (b).

Assertion (a) follows from the fact that for every $x\in
G^{(0)}$ and every compact neighborhood $K$ of $x$, $\pi(K)$ is a
quasi-compact neighborhood of $\pi(x)$.
\pfend


\subsection{Proper actions}
\begin{defi}
Let $G$ be a topological groupoid. Let $Z$ be a topological space
endowed with an action of $G$. Then the action is said to be proper
if $Z\rtimes G$ is a proper groupoid. (We will also say that
$Z$ is a proper $G$-space.)
\end{defi}

A subspace $A$ of a topological space $X$ is said to be relatively
compact (resp. relatively quasi-compact)
if it is included in a compact (resp. quasi-compact)
subspace of $X$. This does not imply that $\overline{A}$ is compact
(resp. quasi-compact).

\begin{prop}\label{prop:proper action}
Let $G$ be a topological groupoid.
Let $Z$ be a topological space endowed with an action
of $G$. Consider the following assertions:
\begin{itemize}
\item[{\rm{(i)}}] $G$ acts properly on $Z$;
\item[{\rm{(ii)}}] $(r,s)\colon Z\rtimes G\to Z\times Z$ is closed and
  $\forall z\in Z$, the stabilizer of $z$ is quasi-compact;
\item[{\rm{(iii)}}] for all quasi-compact subspaces $K$ and $L$ of $Z$,
  $\{g\in G\vert\; Lg\cap K\ne\emptyset\}$ is quasi-compact;
\item[{\rm{(iii)'}}] for all compact subspaces $K$ and $L$ of $Z$,
  $\{g\in G\vert\; Lg\cap K\ne\emptyset\}$ is quasi-compact;
\item[{\rm{(iv)}}] for every quasi-compact subspace $K$ of $Z$, $\{g\in
  G\vert\; Kg\cap K\ne\emptyset\}$ is quasi-compact;
\item[{\rm{(v)}}] there exists a family $(A_i)_{i\in I}$ of subspaces of
  $Z$ such that $Z=\cup_{i\in I}\rond{A}_i$ and $\{g\in G\vert\;A_ig
\cap A_j\ne\emptyset\}$ is relatively quasi-compact for all $i,j\in I$.
\end{itemize}
Then (i)$\iff$(ii)$\implies$(iii)$\implies$(iii)' and (iii)$\implies$(iv).
If $Z$ is locally compact, then (iii)'$\implies$(v) and
(iv)$\implies$(v). If $G^{(0)}$ is Hausdorff and $Z$ is locally
compact Hausdorff, then (i)--(v) are equivalent.
\end{prop}

\pf
(i)$\iff$(ii) follows from Proposition~\ref{prop:proper
  groupoid}[(i)$\iff$(ii)]. Implication
(i)$\implies$(iii) follows from the fact that if $(Z\rtimes
G)_K^L$ is quasi-compact, then its image by the second projection
$Z\rtimes G\to G$ is quasi-compact. (iii)$\implies$(iii)' and
(iii)$\implies$(iv) are obvious.

Suppose that $Z$ is locally compact. Take $A_i\subset Z$ compact such
that $Z=\cup_{i\in I}\rond{A}_i$. If (iii)' is true, then $\{g\in
G\vert\; A_ig\cap A_j\ne\emptyset \}$ is quasi-compact, hence (v). If
(iv) is true, then $\{g\in G\vert\; A_ig\cap A_j\ne\emptyset\}$ is a
subset of the quasi-compact set $\{g\in G\vert\; Kg\cap
K\ne\emptyset\}$, where $K=A_i\cup A_j$, hence (v).

Suppose that $Z$ is locally compact Hausdorff and that
$G^{(0)}$ is Hausdorff. Let us show (v)$\implies$(ii).
Let $C_{ij}$ be a quasi-compact set
such that $\{g\in G\vert\;A_ig\cap A_j\ne\emptyset\}\subset C_{ij}$.

Let $z\in Z$. Choose $i\in I$ such that $z\in A_i$. Since $Z$
and $G^{(0)}$ are Hausdorff,
$\mbox{stab}(z)$ is a closed subspace of $C_{ii}$, therefore it is
quasi-compact.

It remains to prove that the map $\Phi\colon
Z\times_{G^{(0)}}G\to Z\times Z$ given by $\Phi(z,g)=(z,zg)$ is
closed. Let $F\subset Z\times_{G^{(0)}}G$ be a closed subspace, and
$(z,z')\in\overline{\Phi(F)}$. Choose $i$ and $j$ such that $z\in
\rond{A}_i$ and $z'\in \rond{A}_j$. Then
$(z,z')\in\overline{\Phi(F)\cap(A_i\times
  A_j)}\subset\overline{\Phi(F\cap (A_i\times_{G^{(0)}} C_{ij}))}
\subset\overline{\Phi(F\cap(Z\times_{G^{(0)}}C_{ij}))}$.
There exists a net $(z_\lambda,g_\lambda)\in
F\cap (Z\times_{G^{(0)}}C_{ij})$ such
that $(z,z')$ is a limit point of $(z_\lambda,z_\lambda g_\lambda)$.
Since $C_{ij}$ is quasi-compact, after passing to a universal subnet
we may assume that $g_\lambda$ converges to an
element $g\in C_{ij}$. Since $G^{(0)}$ is Hausdorff,
$F\cap (Z\times_{G^{(0)}} C_{ij})$ is
closed in $Z\times C_{ij}$, so $(z,g)$ is an element of $F\cap
(Z\times_{G^{(0)}}C_{ij})$. Using the fact that $Z$ is Hausdorff and
$\Phi$ is continuous, we obtain $(z,z')=\Phi(z,g)\in \Phi(F)$.
\pfend

\begin{rem}
It is possible to define a notion of slice-proper action
which implies properness in the above sense. The two notions
are equivalent in many cases \cite{a74,cem}.
\end{rem}

%

\begin{prop}\label{prop:action of G on G}
Let $G$ be a locally compact groupoid. Then $G$ acts properly on
itself if and only if $G^{(0)}$ is Hausdorff. In particular,
a locally compact space is proper if and only if it is
Hausdorff.
\end{prop}

\pf
It is clear from Proposition~\ref{prop:proper groupoid}(ii) that $G$
acts properly on itself if and only if the product $\varphi\colon
G^{(2)}\to G\times G$ is closed. Since $\varphi$ factors through the
homeomorphism $G^{(2)}\to G\times_{r,r}G$, $(g,h)\mapsto (g,gh)$,
$G$ acts properly on itself if and only if $G\times_{r,r}G$ is
a closed subset of $G\times G$.

If $G^{(0)}$ is Hausdorff, then clearly $G\times_{r,r}G$ is
closed in $G\times G$. Conversely, if $G^{(0)}$ is not Hausdorff, then
there exists $(x,y)\in G^{(0)}\times G^{(0)}$ such that $x\ne y$ and
$(x,y)$ is in the closure of the diagonal of $G^{(0)}\times
G^{(0)}$. It follows that $(x,y)$ is in the closure of
$G\times_{r,r} G$, but $(x,y)\notin  G\times_{r,r} G$,
therefore $G\times_{r,r} G$ is not closed.
\pfend

\subsection{Permanence properties}
\begin{prop}
If $G_1$ and $G_2$ are proper topological groupoids, then $G_1\times
G_2$ is proper.
\end{prop}

\pf
Follows from the fact that the product of two proper maps is proper
\cite[Corollaire I.10.2.3]{bou}.
\pfend

\begin{prop}\label{prop:proper morphism in proper groupoid}
Let $G_1$ and $G_2$ be two topological groupoids such that $G_1^{(0)}$
is Hausdorff and $G_2$ is proper. Suppose that
$f\colon G_1\to G_2$ is a proper morphism. Then $G_1$ is proper.
\end{prop}

\pf
Denote by $r_i$ and $s_i$ the range and source maps of $G_i$
($i=1,2$). Let $\bar f$ be the map $G_1^{(0)}\times G_1^{(0)}\to
G_2^{(0)}\times G_2^{(0)}$ induced from $f$. Since $\bar f\circ (r_1,s_1) =
(r_2,s_2)\circ f$ is proper and $G_1^{(0)}$ is Hausdorff, it follows from
\cite[Proposition I.10.1.5]{bou} that $(r_1,s_1)$ is proper.
\pfend

\begin{prop}\label{prop:proper morphism from proper groupoid}
Let $G_1$ and $G_2$ be two topological groupoids such that $G_1$ is
proper. Suppose that
$f\colon G_1\to G_2$ is a surjective morphism such that the induced
map $f'\colon G_1^{(0)}\to G_2^{(0)}$ is proper. Then $G_2$ is proper.
\end{prop}

\pf
Denote by $r_i$ and $s_i$ the range and source maps of $G_i$
($i=1,2$). Let $F_2\subset G_2$ be a closed subspace, and
$F_1=f^{-1}(F_2)$. Since $G_1$ is proper, $(r_1,s_1)(F_1)$ is closed,
and since $f'\times f'$ is proper, $(f'\times f')\circ (r_1,s_1)(F_1)$
is closed. By surjectivity of $f$, we have $(r_2,s_2)(F_2)=(f'\times
f')\circ (r_1,s_1)(F_1)$. This proves that $(r_2,s_2)$ is
closed. Since for every topological space $T$, the assumptions of the
proposition are also true for the morphism $f\times 1\colon G_1\times
T\to G_2\times T$, the above shows that $(r_2,s_2)\times 1_T$ is
closed. Therefore, $(r_2,s_2)$ is proper.
\pfend

\begin{prop}\label{prop:proper action on product}
Let $G$ be a topological groupoid with $G^{(0)}$ Hausdorff, acting on
two spaces $Y$ and $Z$. Suppose that the action of $G$ on $Z$ is
proper, and that $Y$ is Hausdorff. Then $G$ acts properly on
$Y\times_{G^{(0)}} Z$.
\end{prop}

\pf
The groupoid $(Y\times_{G^{(0)}} Z)\rtimes G$ is isomorphic to the
subgroupoid $\Gamma=\{(y,y',z,g)\in (Y\times Y)\times (Z\rtimes G)
|\; p(y)=r(g),\; y'=yg\}$ of the
proper groupoid $(Y\times Y)\times (Z\rtimes G)$. Since $Y$ and
$G^{(0)}$ are Hausdorff, $\Gamma$ is closed in
$(Y\times Y)\times (Z\rtimes G)$, hence by
Proposition~\ref{prop:proper groupoid}(ii),
$(Y\times_{G^{(0)}} Z)\rtimes G$ is proper.
\pfend

\begin{coro}\label{coro:action of proper groupoid}
Let $G$ be a proper topological groupoid with $G^{(0)}$ Hausdorff.
Then any action of $G$ on a Hausdorff space is proper.
\end{coro}

\pf
Follows from Proposition~\ref{prop:proper action on product} with
$Z=G^{(0)}$.
\pfend

\begin{prop}\label{prop:G and G[T] proper}
Let $G$ be a topological groupoid and $f\colon T\to G^{(0)}$ be
a continuous map. 
\begin{itemize}
\item[{\rm{(a)}}] If $G$ is proper, then $G[T]$ is proper.
\item[{\rm{(ii)}}] If $G[T]$ is proper and $f$ is open surjective,
  then $G$ is proper.
\end{itemize}
\end{prop}

\pf
Let us prove (a).
Suppose first that $T$ is a subspace of $G^{(0)}$ and that $f$ is the
inclusion. Then $G[T]=G_T^T$. Since $(r_T,s_T)$ is the restriction to
$(r,s)^{-1}(T\times T)$ of $(r,s)$, and $(r,s)$ is proper, it follows
that $(r_T,s_T)$ is proper.

In the general case, let $\Gamma=(T\times T)\times G$ and let
$T'\subset T\times G^{(0)}$ be the graph of $f$. Then $\Gamma$ is a
proper groupoid (since it is the product of two proper groupoids), and
$G[T]=\Gamma[T']$.
\par\medskip

Let us prove (b). The only difficulty is to show that $(r,s)$ is
closed. Let $F\subset G$ be a closed subspace and $(y,x)\in
\overline{(r,s)(F)}$. Let $\tilde F= G[T]\cap (T\times T)\times
F$. Choose $(t',t)\in T\times T$ such that $f(t')=y$ and
$f(t)=x$. Denote by $\tilde r$ and $\tilde s$ the range and source
maps of $G[T]$. Then $(t',t)\in\overline{(\tilde r,\tilde s)(\tilde
  F)}$. Indeed, let $\Omega\ni (t',t)$ be an open set, and
$\Omega'=(f\times f)(\Omega)$. Then $\Omega'$ is an open neighborhood
of $(y,x)$, so $\Omega'\cap (r,s)(F)\ne\emptyset$. It follows that
$\Omega\cap (\tilde r,\tilde s)(\tilde F)\ne\emptyset$.

We have proved that $(t',t)\in \overline{(\tilde r,\tilde s)(\tilde
F)}= (\tilde r,\tilde s)(\tilde F)$, so $(y,x)\in (r,s)(F)$.
\pfend

\begin{coro}\label{coro:action on subspace}
Let $G$ be a groupoid acting properly on a topological space $Z$, and
let $Z_1$ be a saturated subspace. Then $G$ acts properly on $Z_1$.
\end{coro}

\pf
Use the fact that $Z_1\rtimes G=(Z\rtimes G)[Z_1]$.
\pfend

\subsection{Invariance by Morita-equivalence}
In this section, we will only consider groupoids whose range maps
are open. We thus need a stability lemma:
%
\begin{lem}\label{lem:stability of r open}
Let $G$ be a topological groupoid whose range map is open. Let $Z$ be
a $G$ space and $f\colon T\to G^{(0)}$ be a continuous open map. Then
the range maps for $Z\rtimes G$ and $G[T]$ are open.
\end{lem}

To prove Lemma~\ref{lem:stability of r open} we need a preliminary
result:

\begin{lem}\label{lem:open maps}
Let $X$, $Y$, $T$ be topological spaces, $g\colon Y\to T$ an open map
and $f\colon X\to T$ continuous. Let $Z=X\times_T Y$.
Then the first projection
$\mbox{pr}_1\colon X\times_T Y\to X$ is open.
\end{lem}

\pf
Let $\Omega\subset Z$ open. There exists an open subspace $\Omega'$ of
$X\times Y$ such that $\Omega=\Omega'\cap Z$. Let $\Delta$ be the
diagonal in $X\times X$. One easily checks that $(\mbox{pr}_1,
\mbox{pr}_1)(\Omega) = (1\times f)^{-1}(1\times g)(\Omega')\cap
\Delta$, therefore $(\mbox{pr}_1, \mbox{pr}_1)(\Omega)$ is open in
$\Delta$. This implies that $\mbox{pr}_1 (\Omega)$ is open in $X$.
\pfend

\begin{proof}[Proof of Lemma~\ref{lem:stability of r open}]
This is clear for $Z\rtimes G=Z\times_{G^{(0)}}G$
using Lemma~\ref{lem:open maps}.

For $G[T]$, first use Lemma~\ref{lem:open maps} to prove that
$T\times_{f,s}G\xrightarrow{pr_2} G$ is open. Since the range map
is open by assumption, the composition
$T\times_{f,s}G\xrightarrow{pr_2} G\xrightarrow{r}G^{(0)}$
is open. Using again Lemma~\ref{lem:open maps},
$G[T]\simeq T\times_{f, r\circ pr_2} (T\times_{f,s}G)
\xrightarrow{pr_1}T$ is open.
\pfend

In order to define the notion of Morita-equivalence for topological
groupoids, we introduce some terminology:

\begin{defi}
Let $G$ be a topological groupoid. Let $T$ be a topological space and
$\rho\colon G^{(0)}\to T$ be a $G$-invariant map.
Then $G$ is said to be $\rho$-proper
if the map $(r,s)\colon G\to G^{(0)}\times_T G^{(0)}$ is
proper. If $G$ acts on a space $Z$ and $\rho\colon Z\to T$ is
$G$-invariant, then the action is said to be $\rho$-proper
if $Z\rtimes G$ is $\rho$-proper.
\end{defi}

It is clear that properness implies $\rho$-properness. There is a
partial converse:

\begin{prop}\label{prop:quasi proper}
Let $G$ be a topological groupoid, $T$ a topological space,
$\rho\colon G^{(0)}\to T$ a $G$-invariant map. If $G$ is $\rho$-proper
and $T$ is Hausdorff, then $G$ is proper.
\end{prop}

\pf
Since $T$ is Hausdorff, $G^{(0)}\times_T G^{(0)}$ is a closed subspace
of $G^{(0)}\times G^{(0)}$, therefore $(r,s)$, being the composition of
the two proper maps $G\to G^{(0)}\times_T G^{(0)}
\to G^{(0)}\times G^{(0)}$, is proper.
\pfend

\begin{rem}
When $T$ is locally Hausdorff, one easily shows that
$G$ is $\rho$-proper iff for every Hausdorff open subspace $V$ of
$T$, $G_{\rho^{-1}(V)}^{\rho^{-1}(V)}$ is proper.
\end{rem}

%
%
%

\begin{prop}\label{prop:def of Morita}\cite{mrw}
Let $G_1$ and $G_2$ be two topological (resp. locally compact)
groupoids. Let $r_i$, $s_i$
($i=1,2$) be the range and source maps of $G_i$, and suppose that
$r_i$ are open. The following are equivalent:
\begin{itemize}
\item[{\rm{(i)}}] there exist a topological (resp. locally compact)
space $T$ and
$f_i\colon T\to G_i^{(0)}$ open surjective
such that $G_1[T]$ and $G_2[T]$ are isomorphic;
\item[{\rm{(ii)}}] there exists a topological (resp. locally compact)
space $Z$, two continuous maps
$\rho\colon Z\to G_1^{(0)}$ and $\sigma\colon Z\to
G_2^{(0)}$, a left action of $G_1$ on $Z$ with momentum map $\rho$ and a
right action of $G_2$ on $Z$ with momentum map $\sigma$ such that
\begin{itemize}
\item[{\rm{(a)}}] the actions commute and are free, the action of
  $G_2$ is $\rho$-proper and the action of $G_1$ is $\sigma$-proper;
\item[{\rm{(b)}}] the natural maps $Z/G_2\to
  G_1^{(0)}$ and $G_1 \backslash Z\to G_2^{(0)}$ induced from $\rho$ and
  $\sigma$ are homeomorphisms.
\end{itemize}
\end{itemize}
Moreover, one may replace (b) by
\begin{itemize}
\item[{\rm (b)'}] $\rho$ and $\sigma$
are open and induce bijections $Z/G_2\to G_1^{(0)}$
and $G_1 \backslash Z\to G_2^{(0)}$.\\
In {\rm (i)}, if $T$ is locally compact then it may be assumed Hausdorff.
\end{itemize}
\end{prop}

If $G_1$ and $G_2$ satisfy the equivalent conditions in
Proposition~\ref{prop:def of Morita}, then they are said to be
Morita-equivalent. Note that if $G_i^{(0)}$ are Hausdorff, then by
Proposition~\ref{prop:quasi proper}, one may replace
``$\rho$-proper'' and ``$\sigma$-proper'' by ``proper''.

To prove Proposition~\ref{prop:def of Morita}, we need preliminary lemmas:
\begin{lem}\label{lem:quotient map open}
Let $G$ be a topological groupoid. The following are equivalent:
\begin{itemize}
\item[{\rm{(i)}}] $r\colon G\to G^{(0)}$ is open;
\item[{\rm{(ii)}}] for every $G$-space $Z$, the canonical mapping
  $\pi\colon Z\to Z/G$ is open.
\end{itemize}
\end{lem}

\pf
To show (ii)$\implies$(i), take $Z=G$: the canonical mapping
$\pi\colon G\to G/G$ is open. Therefore, for every open subspace
$U$ of $G$, $r(U)=G^{(0)}\cap\pi^{-1}(\pi(U))$ is open.

Let us show (i)$\implies$(ii). By Lemma~\ref{lem:stability of r open},
the range map $r\colon Z\rtimes G\to Z$ is open. The conclusion
follows from Proposition~\ref{prop:pi open}.
\pfend

%

\begin{lem}\label{lem:(XxT)/G}
Let $G$ be a topological groupoid such that the range map $r\colon
G\to G^{(0)}$ is open. Let $X$ be a topological space endowed with an
action of $G$ and $T$ a topological space. Then the canonical map
$$f\colon (X\times T)/G\to (X/G)\times T$$
is an isomorphism.
\end{lem}

\pf
Let $\pi\colon X\to X/G$ and $\pi'\colon X\times T\to (X\times T)/G$
be the canonical mappings. Since $\pi$ is open
(Lemma~\ref{lem:quotient map open}),
$f\circ \pi'=\pi\times 1$ is open. Since $\pi'$ is continuous
surjective, it follows that $f$ is open.
\pfend

\begin{lem}\label{lem:bar f proper}
Let $G$ be a topological groupoid whose range map is open and $f\colon
Y\to Z$ a proper, $G$-equivariant map between two $G$-spaces. Then the
induced map $\bar f\colon Y/G\to Z/G$ is proper.
\end{lem}

\pf
We first show that $\bar f$ is closed. Let $\pi\colon Y\to Y/G$ and
$\pi'\colon Z\to Z/G$ be the canonical mappings. Let $A\subset Y/G$ be
a closed subspace. Since $f$ is closed and $\pi$ is continuous,
$(\pi')^{-1}(\bar f(A))=f(\pi^{-1}(A))$ is closed. Therefore, $\bar f
(A)$ is closed.

Applying this to $f\times 1$, we see that for every topological space
$T$, $(Y\times T)/G\to (Z\times T)/G$ is closed. By
Lemma~\ref{lem:(XxT)/G}, $\bar f\times 1_T$ is closed.
\pfend

\begin{lem}\label{lem:composition morphisms}
Let $G_2$ and $G_3$ be topological groupoids whose range maps are open.
Let $Z_1, Z_2$ and $X$ be topological spaces. Suppose there are maps
$$X\xleftarrow{\rho_1} Z_1\xrightarrow{\sigma_1} G_2^{(0)}
\xleftarrow{\rho_2} Z_2\xrightarrow{\sigma_2} G_3^{(0)},$$
a right action of $G_2$ on $Z_1$ with momentum map $\sigma_1$, such that
$\rho_1$ is $G_2$-invariant and the action of $G_2$ is $\rho_1$-proper,
a left action of $G_2$ on $Z_2$ with momentum map
$\rho_2$ and a right $\rho_2$-proper
action of $G_3$ on $Z_2$ with momentum map $\sigma_2$
which commutes with the $G_2$-action.

Then the action of $G_3$ on $Z=Z_1\times_{G_2}Z_2$ is $\rho_1$-proper.
\end{lem}

\begin{proof}
Let $\varphi\colon Z_2\rtimes G_3\to Z_2\times_{G_2^{(0)}} Z_2$
be the map $(z_2,\gamma)\mapsto (z_2,z_2\gamma)$. By
assumption, $\varphi$ is proper, therefore $1_{Z_1}\times \varphi$ is
proper. Let $F=\{(z_1,z_2,z'_2)\in Z_1\times Z_2\times Z_2\vert\;
\sigma_1(z_1)=\rho_2(z_2)=\rho_2(z'_2)\}$. Then $1_{Z_1}\times
\varphi\colon (1\times\varphi)^{-1}(F)\to F$ is proper,
i.e. $Z_1\times_{G_2^{(0)}}(Z_2\rtimes G_3)\to Z_1\times_{G_2^{(0)}}
(Z_2\times_{G_2^{(0)}} Z_2)$ is proper. By
Lemma~\ref{lem:bar f proper}, taking the quotient by $G_2$, we get
that the map
$$\alpha\colon Z\rtimes G_3\to
Z_1\times_{G_2}(Z_2\times_{G_2^{(0)}}Z_2)$$
defined by
$(z_1,z_2,\gamma)\mapsto (z_1,z_2,z_2\gamma)$ is proper.

By assumption, the map $Z_1\rtimes G_2\to Z_1\times_{X}Z_1$
given by $(z_1,g)\mapsto
(z_1,z_1g)$ is proper. Endow $Z_1\rtimes G_2$ with the following right
action of $G_2\times G_2$: $(z_1,g)\cdot(g',g'')=
(z_1g',(g')^{-1}gg'')$. Using
again Lemma~\ref{lem:bar f proper}, the map
\begin{eqnarray*}
\lefteqn{\beta\colon Z_1\times_{G_2}(Z_2\times_{G_2^{(0)}}Z_2)=(Z_1\rtimes
G_2)\times_{G_2\times G_2} (Z_2\times Z_2)}\\
&& \quad \to (Z_1\times_X Z_1)
\times_{G_2\times G_2} (Z_2\times Z_2)\simeq Z\times_X Z
\end{eqnarray*}
is proper. By composition, $\beta\circ\alpha\colon
Z\rtimes G_3\to Z\times_X Z$ is proper.
\pfend

\pf[Proof of Proposition~\ref{prop:def of Morita}]
Let us treat the case of topological groupoids.
Assertion (b') follows from the fact that the canonical mappings
$Z\to Z/G_2$ and $Z\to G_1\backslash Z$ are open
(Lemma~\ref{lem:quotient map open}).

Let us first show that (ii) is an equivalence relation.
Reflexivity is clear (taking $Z=G$, $\rho=r$, $\sigma=s$), and
symmetry is obvious. Suppose that $(Z_1,\rho_1,\sigma_2)$ and
$(Z_2,\rho_2,\sigma_2)$ are equivalences between $G_1$ and $G_2$,
and $G_2$ and $G_3$ respectively. Let $Z=Z_1\times_{G_2}Z_2$
be the quotient of $Z_1\times_{G_2^{(0)}}Z_2$ by the action
$(z_1,z_2)\cdot\gamma = (z_1\gamma,\gamma^{-1}z_2)$ of $G_2$.
Denote by $\rho\colon Z\to G_1^{(0)}$ and $\sigma\colon Z\to
G_3^{(0)}$ the maps induced from $\rho_1\times 1$ and $1\times \sigma_2$.
By Lemma~\ref{lem:open maps}, the first projection
$pr_1\colon Z_1\times_{G_2^{(0)}}Z_2\to Z_1$ is open, therefore
$\rho=\rho_1\circ pr_1$ is open. Similarly, $\sigma$ is open.
It remains to show that the actions of $G_3$ and
$G_1$ are $\rho$-proper and $\sigma$-proper respectively.
For $G_3$, this follows from
Lemma~\ref{lem:composition morphisms}
and the proof for
$G_1$ is similar.

%
This proves that (ii) is an equivalence relation. Now,
let us prove that (i) and (ii) are equivalent.

Suppose (ii). Let $\Gamma=G_1\ltimes Z\rtimes G_2$ and $T=Z$. The maps
$\rho\colon T\to G_1^{(0)}$ and $\sigma\colon T\to G_2^{(0)}$
are open surjective by assumption.
Since $G_1\ltimes Z\simeq Z\times_{G_2^{(0)}} Z$ and
$Z\rtimes G_2\simeq Z\times_{G_1^{(0)}}Z$, we have
$G_2[T]=(T\times T)\times_{G_2^{(0)}\times G_2^{(0)}}G_2
\simeq (Z\rtimes G_2)\times_{s\circ pr_2,\sigma} Z
\simeq (Z\times_{G_1^{(0)}}Z)\times_{\sigma\circ pr_2,\sigma}Z
= Z\times_{G_1^{(0)}} (Z\times_{G_2^{(0)}}Z)
\simeq Z\times_{G_1^{(0)}} (G_1\ltimes Z)
\simeq G_1\ltimes (Z\times_{G_1^{(0)}} Z)
\simeq G_1\ltimes (Z\rtimes G_2)=\Gamma$.
Similarly, $\Gamma\simeq G_1[T]$, hence (i).

Conversely, to prove $(i)\implies (ii)$ it suffices to show that if
$f\colon T\to G^{(0)}$ is open surjective, then $G$ and $G[T]$
are equivalent in the sense (ii), since we know that (ii) is an
equivalence relation. Let $Z=T\times_{r,f}G$.

Let us check that the action of $G$ is $pr_1$-proper. Write
$Z\rtimes G=\{(t,g,h)\in T\times G\times G\vert\; f(t)=r(g)\mbox{ and
  }s(g)=r(h)\}$. One needs to check that the map $Z\rtimes G\to
(T\times_{f,r} G)^2$ defined by $(t,g,h)\mapsto (t,g,t,h)$ is a
homeomorphism onto its image. This follows easily from the facts that
the diagonal map $T\to T\times T$ and the map $G^{(2)}\to G\times G$,
$(g,h)\mapsto (g,gh)$ are homeomorphisms onto their images.

Let us check that the action of $G[T]$ is $s\circ pr_2$-proper.
One easily checks that the groupoid $G'=G[T]\ltimes
(T\times_{f,r}G)$ is isomorphic to a subgroupoid of the trivial
groupoid $(T\times T)\times (G\times G)$. It follows that if $r'$ and $s'$
denote the range and source maps of $G'$, the map $(r',s')$ is a
homeomorphism of $G'$ onto its image.
\par\bigskip
Let us now treat the case of locally compact groupoids.
In the proof that (ii) is a transitive relation,
it just remains to show that $Z$ is locally compact.

Let $U_3$ be a Hausdorff open subspace of
$G_3^{(0)}$. We show that
$\sigma^{-1}(U_3)$ is locally compact. Replacing
$G_3$ by $(G_3)_{U_3}^{U_3}$, we may assume that $G_2$ acts freely and
properly on $Z_2$.
Let $\Gamma$ be the groupoid $(Z_1\times_{G_2^{(0)}} Z_2)\rtimes G_2$,
and $R=(r,s)(\Gamma)\subset (Z_1\times_{G_2^{(0)}}Z_2)^2$. Since the
action of $G_2$ on $Z_2$ is free and proper, there exists a
continuous map $\varphi\colon Z_2\times_{G_3^{(0)}}Z_2 \to G_2$ such
that $z_2=\varphi(z_2,z'_2)z'_2$. Then
$R=\{(z_1,z_2,z'_1,z'_2)\in (Z_1\times_{G_2^{(0)}} Z_2)^2;\;
z'_1=z_1\varphi(z_2,z'_2)\}$ is locally closed. By
Proposition~\ref{prop:G0/G}, $Z=(Z_1\times_{G_2^{(0)}}Z_2)/G$ is
locally compact.

Finally, if (i) holds with $T=\cup_iV_i$ with $V_i$ open Hausdorff,
let $T'=\amalg V_i$. It is clear that $G_1[T']\simeq G_2[T']$.
\pfend

Let us examine standard examples of Morita-equivalences:
\begin{example}\label{ex:G[U]}
Let $G$ be a topological groupoid whose range map is open. Let
$(U_i)_{i\in I}$ be an open cover of $G^{(0)}$ and ${\mathcal{U}}
=\amalg_{i\in I}U_i$. Then $G[{\mathcal{U}}]$ is Morita-equivalent to
$G$.
\end{example}

\begin{example}\label{ex:H1/G/H2}
Let $G$ be a topological groupoid, and let $H_1$, $H_2$ be
subgroupoids such that the range maps $r_i\colon H_i\to H_i^{(0)}$ are
open. Then $(H_1\backslash G^{s(H_1)}_{s(H_2)})\rtimes H_2$ and
$H_1\ltimes (G^{s(H_1)}_{s(H_2)}/H_2)$ are Morita-equivalent.
\end{example}

\pf
Take $Z=G^{s(H_1)}_{s(H_2)}$ and let $\rho\colon Z\to Z/H_2$ and
$\sigma\colon H_1\backslash Z$ be the canonical mappings. The fact that
these maps are open follows from Lemma~\ref{lem:quotient map open}.
\pfend

The following proposition is an immediate consequence of
Proposition~\ref{prop:G and G[T] proper}.

\begin{prop}\label{prop:invariance by Morita}
Let $G$ and $G'$ be two topological groupoids such that the range
maps of $G$ and $G'$ are open. Suppose that $G$ and $G'$
are Morita-equivalent. Then $G$ is proper if and only if $G'$ is proper.
\end{prop}

\begin{coro}
With the notations of Example~\ref{ex:G[U]}, $G$ is proper if and only
if $G[{\mathcal{U}}]$ is proper.
\end{coro}

\section{A topological construction}\label{sec:topological construction}
Let $X$ be a locally compact space. Since $X$ is not necessarily
Hausdorff, a filter\footnote{or a net; we will use indifferently the
two equivalent approaches} $\mathcal{F}$ on $X$ may have more than one
limit. Let $S$ be the set of limits of a convergent filter
$\mathcal{F}$. The goal of this
section is to construct a Hausdorff space ${\mathcal{H}}X$ in which
$X$ is (not continuously) embedded, and such that ${\mathcal{F}}$
converges to $S$ in ${\mathcal{H}}X$.

\subsection{The space ${\mathcal{H}}X$}
\begin{lem}\label{lem:HX}
Let $X$ be a topological space, and $S\subset X$. The following are equivalent:
\begin{itemize}
\item[{\rm{(i)}}] for every family $(V_s)_{s\in S}$ of open sets such
  that $s\in V_s$, and $V_s=X$ except perhaps for finitely many $s$'s,
  one has $\cap_{s\in S} V_s\ne\emptyset$;
\item[{\rm{(ii)}}] for every finite family $(V_i)_{i\in I}$ of open
  sets such that $S\cap V_i\ne\emptyset$ for
  all $i$, one has $\cap_{i\in I}V_i\ne\emptyset$.
\end{itemize}
\end{lem}

\pf
(i)$\implies$(ii): let $(V_i)_{i\in I}$ as in (ii). For all $i$,
choose $s(i)\in S\cap V_i$. Put $W_s=\cap_{s=s(i)}V_i$, with the
convention that an empty intersection is $X$. Then by (i),
$\emptyset\ne\cap_{s\in S}W_s=\cap_{i\in I}V_i$.

(ii)$\implies$(i): let $(V_s)_{s\in S}$ as in (i), and let $I=\{s\in
S\vert\; V_s\ne X\}$. Then $\cap_{s\in S}V_s=\cap_{i\in I}V_i\ne\emptyset$.
\pfend

We shall denote by ${\mathcal{H}}X$ the set of non-empty subspaces $S$ of $X$
which satisfy the equivalent conditions of Lemma~\ref{lem:HX}, and
$\hat{\mathcal{H}}X= {\mathcal{H}}X\cup\{\emptyset\}$.

\begin{lem}\label{lem:S locally finite}
Let $X$ be a locally Hausdorff space. Then every $S\in
{\mathcal{H}}X$ is locally finite. More precisely, if $V$ is a
Hausdorff open subspace of $X$, then $V\cap S$ has at most one element.
\end{lem}

\pf
Suppose $a\ne b$ and $\{a,b\}\subset V\cap S$. Then there exist $V_a$,
$V_b$ open disjoint neighborhoods of $a$ and $b$ respectively; this
contradicts Lemma~\ref{lem:HX}(ii).
\pfend

Suppose that $X$ is locally compact.
We endow $\hat{\mathcal{H}}X$ with a topology. Let us introduce the
notations $\Omega_V=\{S\in {\mathcal{H}}X\vert\; V\cap
S\ne\emptyset\}$ and $\Omega^Q=\{S\in {\mathcal{H}}X\vert\; Q\cap
S=\emptyset\}$. The topology on $\hat{\mathcal{H}}X$ is generated by
the $\Omega_V$'s and $\Omega^Q$'s ($V$ open and $Q$ quasi-compact).
More explicitly,
a set is open if and only if it is a union of sets of the form
$\Omega_{(V_i)_{i\in I}}^Q=\Omega^Q\cap(\cap_{i\in I}\Omega_{V_i})$ where 
$(V_i)_{i\in I}$ is a finite family of open Hausdorff sets and $Q$ is
quasi-compact.

\begin{prop}\label{prop:HX is Hausdorff}
For every locally compact space $X$, the space $\hat{\mathcal{H}}X$ is
Hausdorff.
\end{prop}

\pf
Suppose $S\not\subset S'$ and $S$, $S'\in \hat{\mathcal{H}}X$.
Let $s\in S-S'$. Since $S'$ is locally finite and
since every singleton subspace of $X$ is closed,
there exist $V$ open and $K$ compact such that $s\in V\subset K$ and
$K\cap S'=\emptyset$. Then $\Omega_V$ and
$\Omega^K$ are disjoint neighborhoods of $S$ and $S'$
respectively.
\pfend

For every filter ${\mathcal{F}}$ on $\hat{\mathcal{H}}X$, let
\begin{equation}\label{eqn:L(F)}
L({\mathcal{F}})=\{a\in X\vert\; \forall
V\ni a\mbox{ open}, \Omega_V\in{\mathcal{F}}\}.
\end{equation}

\begin{lem}\label{lem:characterization of limit}
Let $X$ be a locally compact space.
Let ${\mathcal{F}}$ be a filter on $\hat{\mathcal{H}}X$. Then
${\mathcal{F}}$ converges to $S\in\hat{\mathcal{H}}X$ if and only if
properties (a) and (b) below hold:
\begin{itemize}
\item[{\rm{(a)}}] $\forall V$ open, $V\cap S\ne\emptyset\implies
  \Omega_V\in{\mathcal{F}}$;
\item[{\rm{(b)}}] $\forall Q$ quasi-compact, $Q\cap
  S=\emptyset\implies \Omega^Q\in{\mathcal{F}}$.
\end{itemize}
If ${\mathcal{F}}$ is convergent, then $L({\mathcal{F}})$ is its limit.
\end{lem}

\pf
The first statement is obvious, since every open set in
$\hat{\mathcal{H}}X$ is a union of finite intersections of $\Omega_V$'s and
$\Omega^Q$'s.

Let us prove the second statement. It is clear from (a)
that $S\subset L({\mathcal{F}})$.
Conversely, suppose there exists $a\in L({\mathcal{F}})-S$.
Since $S$ is locally finite and every singleton
subspace of $X$ is closed, there exists a
compact neighborhood $K$ of $a$ such that $K\cap S=\emptyset$. Then
$a\in L({\mathcal{F}})$ implies $\Omega_K\in{\mathcal{F}}$, and condition (b)
implies $\Omega^K\in{\mathcal{F}}$, thus
$\emptyset=\Omega^K\cap\Omega_K\in {\mathcal{F}}$, which is
impossible: we have proved the reverse inclusion $L({\mathcal{F}})\subset S$. 
\pfend

\begin{rem}
This means that if $S_\lambda\to S$, then $a\in S$ if and only if
$\forall \lambda$ there exists $s_\lambda\in S_\lambda$ such that
$s_\lambda\to a$.
\end{rem}

\begin{example}
Consider Example~\ref{ex:NH gpd} with $\Gamma=\Z_2$ and $H=\{0\}$.
Then ${\mathcal{H}}G=G\cup \{S\}$ where $S=\{(0,0),(0,1)\}$.
The sequence $(1/n,0)\in G$ converges to $S$ in ${\mathcal{H}}G$,
and $(0,0)$ and $(0,1)$ are two isolated points in ${\mathcal{H}}G$.
\end{example}

\begin{prop}\label{prop:HX is locally compact}
Let $X$ be a locally compact space and $K\subset X$ quasi-compact. Then
$L=\{S\in{\mathcal{H}}X\vert\; S\cap K\ne\emptyset\}$ is compact. The
space ${\mathcal{H}}X$ is locally compact, and it is $\sigma$-compact
if $X$ is $\sigma$-compact.
\end{prop}

\pf
We show that $L$ is compact, and the two remaining assertions follow
easily. Let ${\mathcal{F}}$ be a ultrafilter on $L$. Let
$S_0=L({\mathcal{F}})$. Let us show that
$S_0\cap K\ne\emptyset$: for every $S\in L$, choose
a point $\varphi(S)\in K\cap S$. By quasi-compactness,
$\varphi({\mathcal{F}})$ converges to a point $a\in K$, and it is not
hard to see that $a\in S_0$.

Let us show $S_0\in{\mathcal{H}}X$: let $(V_s)$ ($s\in S_0$) be a
family of open subspaces of $X$ such that $s\in V_s$ for all $s\in S_0$,
and $V_s=X$ for every $s\notin S_1$ ($S_1\subset S_0$ finite). By
definition of $S_0$, $\Omega_{(V_s)_{s\in S_1}}
=\cap_{s\in S_1}\Omega_{V_s}$ belongs to ${\mathcal{F}}$, hence it is
non-empty. Choose
$S\in \Omega_{(V_s)_{s\in S_1}}$, then
$S\cap V_s\ne\emptyset$ for all $s\in
S_1$. By Lemma~\ref{lem:HX}(ii),
$\cap_{s\in S_1}V_s\ne\emptyset$. This shows that $S_0\in{\mathcal{H}}X$.

Now, let us show that ${\mathcal{F}}$ converges to $S_0$.
\begin{itemize}
\item If $V$ is
open Hausdorff such that $S_0\in \Omega_V$, then by definition
$\Omega_V\in {\mathcal{F}}$.
\item If $Q$ is quasi-compact and $S_0\in
\Omega^Q$, then $\Omega^Q\in {\mathcal{F}}$, otherwise one would have
$\{S\in{\mathcal{H}}X\vert\; S\cap Q\ne\emptyset\}\in{\mathcal{F}}$,
which would imply as above that $S_0\cap Q\ne\emptyset$, a
contradiction.
\end{itemize}
From Lemma~\ref{lem:characterization of limit},
${\mathcal{F}}$ converges to $S_0$.
\pfend

\begin{prop}\label{prop:hat HX is compact}
Let $X$ be a locally compact space. Then $\hat{\mathcal{H}}X$ is the
one-point compactification of ${\mathcal{H}}X$.
\end{prop}

\pf
It suffices to prove that $\hat{\mathcal{H}}X$ is compact. The proof
is almost the same as in Proposition~\ref{prop:HX is locally compact}.
\pfend

\begin{rem}
If $f\colon X\to Y$ is a
continuous map from a locally compact space $X$ to any \emph{Hausdorff}
space $Y$, then $f$ induces a continuous map
${\mathcal{H}}f\colon {\mathcal{H}}X\to
Y$. Indeed, for every open subspace $V$ of $Y$,
$({\mathcal{H}}f)^{-1}(V)=\Omega_{f^{-1}(V)}$ is open.
\end{rem}

\begin{prop}\label{prop:extension of actions}
Let $G$ be a topological
groupoid such that $G^{(0)}$ is Hausdorff, and $r\colon G\to
G^{(0)}$ is open.
Let $Z$ be a locally compact space endowed with a continuous action of
$G$. Then ${\mathcal{H}}Z$ is endowed with a continuous action of $G$
which extends the one on $Z$.
\end{prop}

\pf
Let $p\colon Z\to G^{(0)}$ such that $G$ acts on $Z$ with momentum map $p$.
Since $p$ has a continuous extension ${\mathcal{H}}p\colon
{\mathcal{H}}Z\to G^{(0)}$, for all $S\in{\mathcal{H}}Z$,
there exists $x\in G^{(0)}$ such that $S\subset p^{-1}(x)$. For all $g\in
G^x$, write $Sg=\{sg\vert\;s\in S\}$.

Let us show that $Sg\in {\mathcal{H}}Z$. Let $V_s$ ($s\in S$) be open
sets such that $sg\in V_s$. By continuity, there exist open sets
$W_s\ni s$ and $W_g\ni g$ such that for all $(z,h)\in W_s\times_{G^{(0)}}
W_g$, $zh\in V_s$. Let $V'_s=W_s\cap p^{-1}(r(W_g))$. Then $V'_s$ is
an open neighborhood of $s$, so there exists $z\in\cap_{s\in
  S}V'_s$. Since $p(z)\in r(W_g)$, there exists $h\in W_g$ such
that $p(z)=r(h)$. It follows that $zh\in\cap_{s\in S}V_s$. This shows
that $Sg\in {\mathcal{H}}Z$.

Let us show that the action defined above is continuous. Let
$\Phi\colon {\mathcal{H}}Z\times_{G^{(0)}} G\to {\mathcal{H}}Z$ be the action
of $G$ on ${\mathcal{H}}Z$.
Suppose that $(S_\lambda,g_\lambda)\to (S,g)$ and let
$S'=L((S_\lambda,g_\lambda))$. Then for all $a\in S$ there exists
$s_\lambda\in S_\lambda$ such that $s_\lambda\to a$. This implies
$s_\lambda g_\lambda \to ag$, thus $ag\in S'$. The converse
may be proved in a similar fashion, hence $Sg=S'$.

Applying this to any universal net $(S_\lambda,g_\lambda)$
converging to $(S,g)$
and knowing from Proposition~\ref{prop:hat HX is compact} that
$\Phi(S_\lambda,g_\lambda)$ is convergent in $\hat{\mathcal{H}}Z$, we find
that $\Phi(S_\lambda,g_\lambda)$ converges to
$\Phi(S,g)$. This shows that $\Phi$ is continuous in $(S,g)$.

%
%
%
\pfend

\subsection{The space ${\mathcal{H}}'X$}
Let $X$ be a locally compact space.
Let $\Omega'_V=\{S\in {\mathcal{H}}X\vert\; S\subset V\}$. Let
${\mathcal{H}}'X$ be ${\mathcal{H}}X$ as a set, with the coarsest
topology such that the identity map ${\mathcal{H}}'X\to{\mathcal{H}}X$
is continuous, and $\Omega'_V$ is open for every relatively
quasi-compact open set $V$. The space ${\mathcal{H}}'X$
is Hausdorff since ${\mathcal{H}}X$ is Hausdorff, but it is usually not
locally compact.

\begin{lem}\label{lem:card is usc}
Let $X$ be a locally compact space. Then the map
$${\mathcal{H}}'X\to\N^*\cup\{\infty\},\quad
S\mapsto\# S$$
is upper semi-continuous.
\end{lem}

\pf
Let $S\in {\mathcal{H}}'X$ such that $\#S<\infty$. Let $V_s$ ($s\in S$)
be open relatively compact
Hausdorff sets such that $s\in V_s$, and let $W=\cup_{s\in
  S}V_s$. Then $S'\in {\mathcal{H}}'X$ implies $\# (S'\cap V_s)\le
1$, therefore $S'\in\Omega'_W$ implies $\#S'\le\# S$.
\pfend

\begin{prop}\label{prop:property P}
Let $X$ be a locally compact space such that the closure
of every quasi-compact subspace is quasi-compact.
Then
\begin{itemize}
\item[{\rm{(a)}}] the natural map
  ${\mathcal{H}}'X\to{\mathcal{H}}X$ is a homeomorphism,
\item[{\rm{(b)}}] for every compact subspace $K\subset X$, there
  exists $C_K>0$ such that 
$$\forall S\in{\mathcal{H}}X,\;
S\cap K\ne\emptyset\implies \#S\le C_K,$$
\item[{\rm{(c)}}] If $G$ is a locally compact
proper groupoid with $G^{(0)}$ Hausdorff
then $G$ satisfies the above properties.
\end{itemize}
\end{prop}

\pf
%
%
%
To prove (b),
let $K_1$ be a quasi-compact neighborhood of $K$ and let
$K'=\overline{K}_1$. Let $a\in K\cap S$ and suppose there exists
$b\in S- K'$. Then $\rond{K}_1$ and $X-K'$ are disjoint neighborhoods
of $a$ and $b$ respectively, which is impossible.
We deduce that $S\subset K'$.

%
Now, let $(V_i)_{i\in I}$ be a finite cover of $K'$ by open Hausdorff
sets. For all $b\in S$, let $I_b=\{i\in I\vert\; b\in V_i\}$.
By Lemma~\ref{lem:S locally finite}, the
$I_b$'s ($b\in S$) are disjoint, whence one may take $C_K=\# I$.

\par\medskip
To prove (a), denote by $\Delta\subset X\times X$ the diagonal.
Let us first show that $pr_1\colon \overline{\Delta}\to X\times X$
is proper.

Let $K\subset X$ compact. Let $L\subset X$ quasi-compact such that
$K\subset \rond{L}$. If $(a,b)\in \overline{\Delta}\cap (K\times X)$,
then $b\in\overline{L}$: otherwise, $L\times L^c$ would be a
neighborhood of $(a,b)$ whose intersection with $\Delta$ is
empty. Therefore, $pr_1^{-1}(K)=\overline{\Delta}\cap (K\times
\overline{L})$ is quasi-compact, which shows that $pr_1$ is proper.

%
It remains to prove that $\Omega'_V$
is open in ${\mathcal{H}}X$ for every relatively quasi-compact
open set $V\subset X$. Let $S\in
\Omega'_V$, $a\in S$ and $K$ a compact neighborhood of $a$. Let
$L=pr_2(\overline{\Delta}\cap (K\times X))$. Then $Q=L-V$ is
quasi-compact, and $S\in \Omega^Q_{\rond{K}}\subset\Omega'_V$,
therefore $\Omega'_V$ is a neighborhood of each of its points.
%
%
%
%
\par\medskip
To prove (c),
let $K\subset G$ be a quasi-compact subspace.
Then $L={r(K)\cup s(K)}$ is quasi-compact,
thus $G_L^L$ is also quasi-compact. But
$\overline{K}$ is closed and $\overline{K}\subset G_L^L$,
therefore $\overline{K}$ is quasi-compact.
\pfend

\section{Haar systems}\label{sec:haar systems}
\subsection{The space $C_c(X)$}
For every locally compact space $X$,
$C_c(X)_0$ will denote the set of functions $f\in C_c(V)$ ($V$ open
Hausdorff), extended by 0 outside $V$. Let
$C_c(X)$ be the linear span of $C_c(X)_0$.
Note that functions in $C_c(X)$ are not necessarily continuous.

\begin{prop}\label{prop:C_c(X)}
Let $X$ be a locally compact space, and let
$f\colon X\to \C$. The following are equivalent:
\begin{itemize}
\item[{\rm{(i)}}] $f\in C_c(X)$;
\item[{\rm{(ii)}}] $f^{-1}(\C^*)$ is relatively quasi-compact,
and for every filter ${\mathcal{F}}$ on $X$, let
$\tilde{\mathcal{F}}=i({\mathcal{F}})$, where $i\colon X\to
{\mathcal{H}}X$ is the canonical inclusion; if $\tilde{\mathcal{F}}$
converges to $S\in {\mathcal{H}}X$, then $\lim_{\mathcal{F}} f =
\sum_{s\in S}f(s)$.
\end{itemize}
\end{prop}

\pf
Let us show (i)$\implies$(ii). By linearity, it is enough to consider
the case $f\in C_c(V)$, where $V\subset X$ is open Hausdorff. Let
$K$ be the compact set $\overline{f^{-1}(\C^*)}\cap V$. Then
$f^{-1}(\C^*)\subset K$. Let ${\mathcal{F}}$ and $S$ as in (ii). If
$S\cap V=\emptyset$, then $S\in \Omega^K$, hence $\Omega^K\in
\tilde{\mathcal{F}}$, i.e. $X-K\in {\mathcal{F}}$. Therefore,
$\lim_{\mathcal{F}} f=0=\sum_{s\in S}f(s)$.

If $S\cap V=\{a\}$, then $a$ is a limit point of ${\mathcal{F}}$,
therefore $\lim_{\mathcal{F}}f=f(a)=\sum_{s\in S}f(s)$.

Let us show (ii)$\implies$(i) by induction on $n\in \N^*$ such that
there exist $V_1,\ldots V_n$ open Hausdorff and $K$ quasi-compact
satisfying $f^{-1}(\C^*)\subset K\subset V_1\cup\cdots\cup V_n$.

For $n=1$, for every $x\in V_1$, let ${\mathcal{F}}$ be a ultrafilter
convergent to $x$. By Proposition~\ref{prop:hat HX is compact},
$\tilde{\mathcal{F}}$ is convergent; let $S$ be its limit, then
$\lim_{\mathcal{F}}f = \sum_{s\in S}f(s)=f(x)$, thus $f_{\vert V_1}$
is continuous.

Now assume the implication is true for $n-1$ ($n\ge 2$) and let us
prove it for $n$.
Since $K$ is quasi-compact, there
exist $V'_1,\ldots,V'_n$ open sets, $K_1\ldots,K_n$ compact such that
$K\subset V'_1\cup\cdots\cup V'_n$ and $V'_i\subset K_i\subset
V_i$.
Let $F=(V'_1\cup\cdots\cup V'_n)-(V'_1\cup\cdots\cup
V'_{n-1})$. Then $F$ is closed in $V'_n$ and
$f_{\vert F}$ is continuous. Moreover,
$f_{\vert F}=0$ outside $K'=K-(V'_1\cup\cdots\cup V'_{n-1})$ which
is closed in $K$, hence quasi-compact, and Hausdorff, since $K'\subset
V'_n$. Therefore, $f_{\vert F}\in C_c(F)$. It follows that there exists
an extension $h\in C_c(V'_n)$ of $f_{\vert F}$. By considering $f-h$,
we may assume that $f=0$ on $F$, so $f=0$ outside $K'=K_1\cup\cdots\cup
K_{n-1}$. But $K'\subset V_1\cup\cdots\cup V_{n-1}$, hence by
induction hypothesis, $f\in C_c(X)$.
\pfend

\begin{coro}\label{cor:limit in C_c}
Let $X$ be a locally compact space, $f\colon X\to \C$, $f_n\in C_c(X)$.
Suppose that there exists fixed quasi-compact set $Q\subset X$
such that $f_n^{-1}(\C^*)\subset Q$ for all $n$, and
$f_n$ converges uniformly to $f$. Then $f\in C_c(X)$.
\end{coro}

\begin{lem}\label{lem:predetermined cover}
Let $X$ be a locally compact space. Let $(U_i)_{i\in I}$ be an open cover of
$X$ by Hausdorff subspaces. Then every $f\in C_c(X)$ is a finite sum
$f=\sum f_i$, where $f_i\in C_c(U_i)$.
\end{lem}

\pf
See \cite[Lemma 1.3]{ks02}.
\pfend

\begin{lem}\label{lem:density C_c(X)xC_c(Y)}
Let $X$ and $Y$ be locally compact spaces. Let $f\in C_c(X\times
Y)$. Let $V$ and $W$ be open subspaces of $X$ and $Y$ such that
$f^{-1}(\C^*)\subset Q \subset V\times W$ for some quasi-compact
set $Q$. Then there exists a sequence
$f_n\in C_c(V)\otimes C_c(W)$ such that $\lim_{n\to\infty}
\|f-f_n\|_\infty=0$.
\end{lem}

\pf
We may assume that $X=V$ and $Y=W$.
Let $(U_i)$ (resp. $(V_j)$) be an open cover of $X$ (resp. $Y$) by
Hausdorff subspaces. Then every element of $C_c(X\times Y)$ is a linear
combination of elements of $C_c(U_i\times V_j)$
(Lemma~\ref{lem:predetermined cover}). The conclusion follows from the
fact that the image of $C_c(U_i)\otimes C_c(V_j)\to C_c(U_i\times
V_j)$ is dense.
\pfend

\begin{lem}\label{lem:C_c(X) to C_c(Y)}
Let $X$ be a locally compact space and $Y\subset X$ a closed
subspace. Then the restriction map $C_c(X)\to C_c(Y)$ is well-defined
and surjective.
\end{lem}

\pf
Let $(U_i)_{i\in I}$ be a cover of $X$ by Hausdorff open subspaces. The map
$C_c(U_i)\to C_c(U_i\cap Y)$ is surjective (since $Y$ is closed), and
$\oplus_{i\in I} C_c(U_i\cap Y)\to C_c(Y)$ is surjective
(Lemma~\ref{lem:predetermined cover}). Therefore, the map
$\oplus_{i\in I} C_c(U_i)\to C_c(Y)$ is surjective. Since
it is also the composition of the surjective map
$\oplus_{i\in I} C_c(U_i)\to
C_c(X)$ and of the restriction map $C_c(X)\to C_c(Y)$, the conclusion
follows.
\pfend

%
%
%

\subsection{Haar systems}
Let $G$ be a locally compact proper groupoid with Haar system (see
definition below) such that $G^{(0)}$ is Hausdorff. If $G$ is Hausdorff, then
$C_c(G^{(0)})$ is endowed with the $C^*_r(G)$-valued scalar product
$\langle\xi,\eta\rangle (g)= \overline{\xi(r(g))}\eta(s(g))$. Its
completion is a $C^*_r(G)$-Hilbert module. However, if $G$ is not
Hausdorff, the function $g\mapsto \overline{\xi(r(g))}\eta(s(g))$ does
not necessarily belong to $C_c(G)$, therefore we need a different
construction in order to obtain a $C^*_r(G)$-module.

\begin{defi}\cite[pp. 16-17]{ren}
Let $G$ be a locally compact groupoid such that $G^x$ is Hausdorff for
every $x\in G^{(0)}$. A Haar system is a
family of positive measures $\lambda=\{\lambda^x\vert\; x\in
G^{(0)}\}$ such that $\forall x,y\in G^{(0)}$, $\forall\varphi\in C_c(G)$,
\begin{itemize}
\item[{\rm{(i)}}] ${\mathrm {supp}}(\lambda^x)=G^x$;
\item[{\rm{(ii)}}]
$\lambda(\varphi)\colon x\mapsto \int_{g\in
  G^x}\varphi(g)\,\lambda^x({\mathrm{d}}g)\quad\in C_c(G^{(0)});$
\item[{\rm{(iii)}}]
$\int_{h\in G^x}\varphi(gh)\,\lambda^x({\mathrm{d}}h) =
\int_{h\in G^y} \varphi(h)\,\lambda^y({\mathrm{d}}h)$.
\end{itemize}
\end{defi}

Note that $G^x$ is automatically Hausdorff if
$G^{(0)}$ is Hausdorff (Prop.~\ref{prop:Gx Hausdorff}). Recall
also \cite[p. 36]{pat} that the range map for $G$ is open.

\begin{lem}\label{lem:measure of quasi-compact finite}
Let $G$ be a locally compact groupoid with Haar system. Then for every
quasi-compact subspace $K$ of $G$, $\sup_{x\in G^{(0)}}\lambda^x (K\cap
G^x)<\infty$.
\end{lem}

\pf
It is easy to show that there exists $f\in C_c(G)$ such that $1_K\le
f$. Since $\sup_{x\in G^{(0)}}\lambda(f)(x)<\infty$,
the conclusion follows.
\pfend

\begin{lem}\label{lem:Haar product}
Let $G$ be a locally compact groupoid with Haar system such that
$G^{(0)}$ is Hausdorff. Suppose that $Z$ is a locally compact space
and that $p\colon Z\to G^{(0)}$ is continuous. Then for every $f\in
C_c(Z\times_{p,r}G)$, $\lambda(f)\colon
z\mapsto\int_{g\in G^{p(z)}} f(z,g)
\,\lambda^{p(z)}({\mathrm{d}}g)$ belongs to $C_c(Z)$.
\end{lem}

\pf
By Lemma~\ref{lem:C_c(X) to C_c(Y)}, $f$ is the restriction
of an element of $C_c(Z\times G)$.

If $f(z,g)=f_1(z)f_2(g)$, then $\psi(x)=\int_{g\in G^x}
f_2(g)\,\lambda^x(dg)$ belongs to $C_c(G^{(0)})$, therefore $\psi\circ
p\in C_b(Z)$. It follows that $\lambda(f)=f_1(\psi\circ p)$ belongs to
$C_c(Z)$.

By linearity, if $f\in C_c(Z)\otimes C_c(G)$, then $\lambda(f)\in C_c(Z)$.

Now, for every $f\in C_c(Z\times G)$, there exist relatively
quasi-compact open subspaces $V$ and $W$ of $Z$ and $G$ and a sequence
$f_n\in C_c(V)\otimes C_c(W)$ such that $f_n$ converges uniformly to
$f$. From Lemma~\ref{lem:measure of quasi-compact finite},
$\lambda(f_n)$ converges uniformly to $\lambda(f)$, and
$\lambda(f_n)\in C_c(Z)$. From Corollary~\ref{cor:limit in C_c},
$\lambda(f)\in C_c(Z)$.
\pfend

\begin{prop}
Let $G$ be a locally compact groupoid with Haar system such that
$G^{(0)}$ is Hausdorff. If $G$ acts on a locally compact space $Z$
with momentum map $p\colon Z\to G^{(0)}$, then $(\lambda^{p(z)})_{z\in Z}$ is a
Haar system on $Z\rtimes G$.
\end{prop}

\pf
Results immediately from Lemma~\ref{lem:Haar product}.
\pfend

\section{The Hilbert module of a proper groupoid}
\label{sec:Hilbert module of proper groupoid}
\subsection{The space $X'$}
Before we construct a Hilbert module associated to a proper groupoid,
we need some preliminaries.
Let $G$ be a locally compact groupoid
such that $G^{(0)}$ is Hausdorff. Denote
by $X'$ the closure of $G^{(0)}$ in ${\mathcal{H}}G$.

\begin{lem}\label{lem:S is a group}
Let $G$ be a locally compact groupoid
such that $G^{(0)}$ is Hausdorff.
Then for all $S\in X'$, $S$ is a subgroup of $G$.
\end{lem}

\pf
Since $r$ and $s\colon G\to G^{(0)}$ extend continuously to maps
${\mathcal{H}}G\to G^{(0)}$, and since $r=s$ on $G^{(0)}$, one has
${\mathcal{H}}r={\mathcal{H}}s$ on $X'$,
i.e. $\exists x_0\in G^{(0)}$, $S\subset G_{x_0}^{x_0}$.

Let ${\mathcal{F}}$ be a filter on $G^{(0)}$ whose limit is $S$. Then
$a\in S$ if and only if $a$ is a limit point of ${\mathcal{F}}$. Since
for every $x\in G^{(0)}$ we have $x^{-1}x=x$, it follows that for every $a$,
$b\in S$ one has $a^{-1}b\in S$, whence $S$ is a subgroup of
$G_{x_0}^{x_0}$.
\pfend

Denote by $q\colon X'\to G^{(0)}$ the map such that $S\subset
G_{q(S)}^{q(S)}$. The map $q$ is continuous since it is the
restriction to $X'$ of ${\mathcal{H}}r$.

\begin{lem}\label{lem:S/S_0}
Let $G$ be a locally compact proper groupoid
such that $G^{(0)}$ is Hausdorff. Let ${\mathcal{F}}$ be a filter on
$X'$, convergent to $S$. Suppose that
$q({\mathcal{F}})$ converges to $S_0\in X'$. Then
$S_0$ is a normal subgroup of $S$, and there exists $\Omega\in
{\mathcal{F}}$ such that $\forall S'\in\Omega$, $S'$ is
group-isomorphic to $S/S_0$. In particular, $\{S'\in X'\vert\;
\#S=\#S_0\# S'\}\in {\mathcal{F}}$.
\end{lem}

\pf
Using Proposition~\ref{prop:property P}, we see that $S$ is finite.

We shall use the notation $\tilde\Omega_{(V_i)_{i\in
    I}}=\Omega_{(V_i)_{i\in I}}\cap\Omega'_{\cup_{i\in I}V_i}$.
Let $V'_s\subset V_s$ ($s\in S$) be Hausdorff, open neighborhoods of
$s$, chosen small enough so that for some $\Omega\in{\mathcal{F}}$,
\begin{itemize}
\item[{\rm{(a)}}] $\Omega\subset\tilde\Omega_{(V'_s)_{s\in S}}$;
\item[{\rm{(b)}}] $V'_{s_1}V'_{s_2}\subset V_{s_1s_2}$, $\forall s_1$,
  $s_2\in S$.
\item[{\rm{(c)}}] $\forall s\in S-S_0$, $\forall S'\in\Omega$,
  $q(S')\notin V_s$;
\item[{\rm{(d)}}] $q(\Omega)\subset \tilde\Omega_{(V_s)_{s\in S_0}}$;

\end{itemize}
Let $S'\in \Omega$. Let $\varphi\colon S\to S'$ such that
$\{\varphi(s)\}=S'\cap V'_s$. Then $\varphi$ is well-defined since
$S'\cap V'_s\ne\emptyset$ (see (a)) and $V'_s$ is Hausdorff.

If $s_1$, $s_2\in S$ then $\varphi(s_i)\in S'\cap V'_{s_i}$. By (b),
$\varphi(s_1)\varphi(s_2)\in S'\cap V_{s_1s_2}$. Since $V_{s_1s_2}$ is
Hausdorff and also contains $\varphi(s_1s_2)\in S'$, we have
$\varphi(s_1s_2)= \varphi(s_1)\varphi(s_2)$. This shows that $\varphi$
is a group morphism.

The map $\varphi$ is surjective, since $S'\subset\cup_{s\in S} V'_s$
(see (a)).

By (c), $\ker(\varphi)\subset S_0$ and by (d), $S_0\subset\ker(\varphi)$.
\pfend

%

Suppose now that the range map $r\colon G\to G^{(0)}$ is open.
Then $X'$ is endowed with an action of $G$
(Prop.~\ref{prop:extension of actions})
defined by $S\cdot g= g^{-1}Sg=\{g^{-1}sg\vert\; s\in S\}$.
%
%
%
%
%

\subsection{Construction of the Hilbert module}
Now, let $G$ be a locally compact, proper groupoid. Assume that
$G$ is endowed with a Haar system, and that $G^{(0)}$ is
Hausdorff. Let
$${\mathcal{E}}^0=\{f\in C_c(X')\vert\;
f(S)=\sqrt{\#S}f(q(S))\;\forall S\in X'\}.$$
($q(S)\in G^{(0)}$ is identified to $\{q(S)\}\in X'$.)

Define, for all $\xi$, $\eta\in {\mathcal{E}}^0$ and $f\in C_c(G)$:
$\langle\xi,\eta\rangle (g)=\overline{\xi(r(g))}\eta(s(g))$ and
%
%
%
$$(\xi f)(S)=\int_{g\in G^{q(S)}} \xi(g^{-1}Sg) f(g^{-1})\,\lambda^x(dg).$$
\begin{prop}
With the above assumptions,
the completion ${\mathcal{E}}(G)$ of
${\mathcal{E}}^0$ with respect to the norm
$\|\xi\|=\|\langle\xi,\xi\rangle\|^{1/2}$ is a $C^*_r(G)$-Hilbert module.
\end{prop}

We won't give the direct proof here since this is a particular
case of Theorem~\ref{thm:correspondence}
(see Example~\ref{ex:proper-morphism}(c)).

\section{Cutoff functions}\label{sec:cutoff}
If $G$ is a locally compact Hausdorff proper groupoid with Haar
system. Assume for simplicity that $G^{(0)}/G$ is compact. Then there
exists a so-called ``cutoff'' function $c\in C_c(G^{(0)})_+$ such that
for every $x\in G^{(0)}$, $\int_{g\in G^x}
c(s(g))\,\lambda^x({\mathrm{d}}g)=1$, and the function $g\mapsto
\sqrt{c(r(g))c(s(g))}$ defines projection in $C^*_r(G)$. However, if $G$
is not Hausdorff, then the above function does not belong to $C_c(G)$
is general, thus we need another definition of a cutoff function.
\par\medskip
Let $X'_{\ge k}=\{S\in X'\vert\;\# S\ge k\}$.
By Lemma~\ref{lem:card is usc}, $X'_{\ge k}$ is closed.

\begin{lem}
Let $G$ be a locally compact, proper groupoid with $G^{(0)}$
Hausdorff. Let $X_{\ge k}=q(X'_{\ge k})$.
Then $X_{\ge k}$ is closed in $G^{(0)}$.
\end{lem}

\pf
It suffices to show that for every compact subspace $K$ of $G^{(0)}$, $X_{\ge
  k}\cap K$ is closed. Let $K'=G_K^K$. Then $K'$ is quasi-compact, and
from Proposition~\ref{prop:HX is locally compact}, $K''=\{S\in
{\mathcal{H}}G\vert\; S\cap K'\ne\emptyset\}$ is compact. The set
$q^{-1}(K)\cap X'_{\ge k}=K''\cap X'_{\ge k}$ is closed in $K''$,
hence compact; its image by $q$ is $X_{\ge k}\cap K$.
\pfend

\begin{lem}\label{lem:local construction of c}
Let $G$ be a locally compact, proper groupoid, with $G^{(0)}$
Hausdorff. Let $\alpha\in\R$.
For every compact set $K\subset G^{(0)}$, there exists $f\colon
X'_K\to \R_+^*$ continuous, where $X'_K=q^{-1}(K)\subset X'$, such
that
$$\forall S\in X'_K,\quad f(S)=f(q(S))(\# S)^\alpha.$$
\end{lem}

\pf
Let $K'=G_K^K$. It is closed and quasi-compact. From
Proposition~\ref{prop:HX is locally compact}, $X'_K$ is quasi-compact.
For every $S\in X'_K$,
we have $S\subset K'$. By Proposition~\ref{prop:property P}, there
exists $n\in \N^*$ such that $X'_{\ge n+1}\cap X'_K=\emptyset$. We can
thus proceed by reverse induction: suppose constructed $f_{k+1}\colon
X'_K\cap q^{-1}(X_{\ge k+1}) \to \R_+^*$ continuous such that
$f_{k+1}(S)=f_{k+1}(q(S))(\# S)^\alpha$ for all $S\in
X'_K\cap q^{-1}(X_{\ge k+1})$.

Since $X'_K\cap q^{-1}(X_{\ge k+1})$ is closed in the compact set
$X'_K\cap q^{-1}(X_{\ge k})$, there exists a continuous extension
$h\colon X'_K\cap q^{-1}(X_{\ge k})\to \R$ of $f_{k+1}$.
Replacing $h(x)$ by $\sup(h(x),\inf f_{k+1})$, we may assume that
$h(X'_K\cap q^{-1}(X_{\ge k}))\subset \R_+^*$. Put
$f_k(S)=h(q(S))(\# S)^\alpha$. Let us show that $f_k$ is continuous.

Let ${\mathcal{F}}$ be a ultrafilter on $X'_K\cap q^{-1}(X_{\ge
  k})$, and let $S$ be its limit. Since $q({\mathcal{F}})$ is a
  ultrafilter on $K$, it has a limit $S_0\in X'_K$.

For every $S_1\in q^{-1}(X_{\ge k})$, choose $\psi(S_1)\in X'_{\ge k}$
such that $q(S_1)=q(\psi(S_1))$. Let $S'\in X'_K\cap X'_{\ge k}$
be the limit of $\psi({\mathcal{F}})$.

From Lemma~\ref{lem:S/S_0},
$\Omega_1=\{ S_1\in X'_K\cap q^{-1}(X_{\ge k})\vert\; \#S = \#
S_0\#S_1\}$ is an element of ${\mathcal{F}}$, and $\Omega_2=\{
S_2\in X'_{\ge k}\vert\; \#S'=\#S_0\#S_2\}$ is an element of
$\psi({\mathcal{F}})$.

\begin{itemize}
\item If $\# S_0>1$, then $S'\in X_{\ge k+1}$, so $S$
  and $S_0$ belong to $q^{-1}(X_{\ge k+1})$. Therefore,
$f_k(S_1)=(\#S_1)^\alpha h(q(S_1))$ converges with respect to
${{\mathcal{F}}}$ to
\begin{eqnarray*}
\lefteqn{
\frac{(\#S)^\alpha}{(\#S_0)^\alpha} h(S_0)=
\frac{(\#S)^\alpha}{(\#S_0)^\alpha}
f_{k+1}(S_0)=f_{k+1}(S)}\\
&&=f_{k+1}(q(S)) (\#S)^\alpha
= h(q(S)) (\#S)^\alpha=f_k(S).
\end{eqnarray*}
\item If $S_0=\{q(S)\}$, then $f_k(S_1)=(\#S_1)^\alpha
h(q(S_1))$ converges with respect to ${{\mathcal{F}}}$ to
$(\#S)^\alpha h(q(S))= f_k(S)$.
\end{itemize}
Therefore, $f_k$ is a continuous extension of $f_{k+1}$.
\pfend

\begin{theo}\label{thm:cutoff}
Let $G$ be a locally compact, proper groupoid such
that $G^{(0)}$ is Hausdorff and $G^{(0)}/G$ is $\sigma$-compact.
Let $\pi\colon G^{(0)}\to G^{(0)}/G$ be the canonical mapping.
Then there exists $c\colon X'\to \R_+$ continuous such that
\begin{itemize}
\item[{\rm{(a)}}] $c(S)=c(q(S))\#S$ for all $S\in X'$;
\item[{\rm{(b)}}] $\forall \alpha\in G^{(0)}/G$, $\exists x\in
  \pi^{-1}(\alpha)$, $c(x)\ne 0$;
\item[{\rm{(c)}}] $\forall K\subset G^{(0)}$ compact, $\mathrm{supp}(c)\cap
    q^{-1}(F)$ is compact, where $F=s(G^K)$.
\end{itemize}
If moreover $G$ admits a Haar system, then there exists $c\colon X'
\to \R_+$ continuous satisfying (a), (b), (c) and
\begin{itemize}
\item[{\rm{(d)}}]
$\displaystyle
\forall x\in G^{(0)},\quad \int_{g\in G^x} c(s(g))\,\lambda^x(dg)=1$.
\end{itemize}
\end{theo}

\pf
There exists a locally
finite cover $(V_i)$ of $G^{(0)}/G$ by relatively compact open
subspaces. Since $\pi$ is open and $G^{(0)}$ is locally compact, there exists
$K_i\subset G^{(0)}$ compact such that $\pi(K_i)\supset V_i$. Let
$(\varphi_i)$ be a partition of unity associated to the cover
$(V_i)$. For every $i$, from Lemma~\ref{lem:local construction of c},
there exists $c_i\colon X'_{K_i}\to \R_+^*$ continuous such that
$c_i(S)=c_i(q(S))\#S$ for all $S\in X'_{K_i}$. Let
$$c(S)=\sum_i c_i(S)\varphi_i(\pi(q(S))).$$
It is clear that $c$ is continuous from $X'$ to $\R_+$, and that
$c(S)=c(q(S))\# S$.

Let us prove (b): let $x_0\in G^{(0)}$. There exists $i$ such that
$\varphi_i(\pi(x_0))\ne 0$. Choose $x\in K_i$
such that $\pi(x)=\pi(x_0)$, then $c(x)\ge
c_i(x)\varphi_i(\pi(x_0))>0$.

Let us show (c). Note that $F=\pi^{-1}(\pi(K))$ is closed, so
$q^{-1}(F)$ is closed.  Let $K_1$ be a compact neighborhood of $K$
and $F_1=\pi^{-1}(\pi(K_1))$.
Let $J=\{ i\vert\; V_i\cap \pi(K_1)\ne\emptyset\}$.
Then for all $i\notin J$,
$c_i(\varphi_i\circ \pi\circ q)=0$ on $q^{-1}(F_1)$, therefore
$c=\sum_{j\in J} c_j(\varphi_j\circ\pi\circ q)$ in a neighborhood of
$q^{-1}(F)$.
Since for all $i$, ${\mathrm{supp}}(c_i(\varphi_i\circ
\pi\circ q))$ is compact and since $J$ is finite,
${\mathrm{supp}}(c)\cap
q^{-1}(F)\subset \cup_{i\in J}{\mathrm{supp}}(c_i(\varphi_i\circ
\pi\circ q))$ is compact.

Let us show the last assertion. Let $\varphi(g)=c(s(g))$. Let
${\mathcal{F}}$ be a filter on $G$ convergent in ${\mathcal{H}}G$ to
$A\subset G$. Choose $a\in A$ and let $S=a^{-1}A$. Then
$s({\mathcal{F}})$ converges to $S$ in ${\mathcal{H}}G$, hence
$$\lim_{{\mathcal{F}}}\varphi=\#S c(s(a))=\sum_{g\in S}
c(s(g))=\sum_{g\in S}\varphi(g).$$
For every compact set $K\subset G^{(0)}$,
\begin{eqnarray*}
\lefteqn{\{g\in G\vert\; r(g)\in K\mbox{ and }\varphi(g)\ne 0\}}\\
&\subset& \{g\in G\vert\; r(g)\in K\mbox{ and
  }s(g)\in{\mathrm{supp}}(c)\}\\
&\subset& G^K_{q({\mathrm{supp}}(c)\cap q^{-1}(F))},
\end{eqnarray*}
so $G^K\cap \{g\in G\vert\;\varphi(g)\ne 0\}$ is included in a
quasi-compact set. Therefore, for every $l\in C_c(G^{(0)})$, $g\mapsto
l(r(g))\varphi(g)$ belongs to $C_c(G)$. It follows that
$h(x)=\int_{g\in G^x}\varphi(g)\,\lambda^x(dg)$ is a continuous
function. Moreover, for every $x\in G^{(0)}$ there exists $g\in G^x$
such that $\varphi(g)\ne 0$, so $h(x)>0$ $\forall x\in G^{(0)}$. It
thus suffices to replace $c(x)$ by $c(x)/h(x)$.
\pfend

\begin{example}
In Example~\ref{ex:NH gpd} with $\Gamma=\Z_n$ and $H=\{0\}$,
the cutoff function is the unique continuous extension to $X'$
of the function $c(x)=1$ for $x\in (0,1]$, and
$c(0)=1/n$.
\end{example}

\begin{prop}\label{prop:E_G is projective}
Let $G$ be a locally compact, proper groupoid with
Haar system such that $G^{(0)}$ is Hausdorff and $G^{(0)}/G$ is
compact. Let $c$ be a cutoff function.
Then the function $p(g)=\sqrt{c(r(g))c(s(g))}$ defines a selfadjoint
projection $p\in C^*_r(G)$, and ${\mathcal{E}}(G)$ is isomorphic to
$pC^*_r(G)$.
\end{prop}

\pf
Let $\xi_0(x)=\sqrt{c(x)}$. Then one easily checks that
$\xi_0\in {\mathcal{E}}^0$, $\langle\xi_0,\xi_0\rangle = p$ and
$\xi_0\langle\xi_0,\xi_0\rangle =\xi_0$, therefore $p$ is a
selfadjoint projection in $C^*_r(G)$. The maps
\begin{eqnarray*}
{\mathcal{E}}(G)\rightarrow pC^*_r(G),&\quad&
\xi\mapsto \langle\xi_0,\xi\rangle = p \langle\xi_0,\xi\rangle\\
pC^*_r(G)\rightarrow {\mathcal{E}}(G),&\quad&
a\mapsto \xi_0a=\xi_0pa
\end{eqnarray*}
are inverses from each other.
\pfend

\section{Generalized morphisms and $C^*$-algebra correspondences}
\label{sec:correspondences}
%
{\bf Until the end of the paper, all groupoids are assumed locally compact,
with open range map.}
In this section, we introduce a notion of generalized morphism
for locally compact groupoids which are not necessarily Hausdorff,
and a notion of locally proper generalized morphism.

Then, we show that a locally proper generalized
morphism from $G_1$ to $G_2$ which satisfies an
additional condition induces a $C^*_r(G_1)$-module ${\mathcal{E}}$
and a $*$-morphism $C^*_r(G_2)\to {\mathcal{K}}({\mathcal{E}})$,
hence an element of $KK(C^*_r(G_2),C^*_r(G_1))$.

\subsection{Generalized morphisms}
\begin{defi}\label{def:generalized morphism}\cite{Hae,hs,lan2,leg,moe,mrw}
Let $G_1$ and $G_2$ be two groupoids. A generalized
morphism from $G_1$ to $G_2$ is a triple $(Z,\rho,\sigma)$ where
$$G_1^{(0)}\xleftarrow{\rho}Z\xrightarrow{\sigma} G_2^{(0)},$$
$Z$ is endowed with a left action of $G_1$ with momentum map
$\rho$ and a right action of
$G_2$ with momentum map $\sigma$ which commute, such that
\begin{itemize}
\item[{\rm{(a)}}] the action of $G_2$ is free and $\rho$-proper,
\item[{\rm{(b)}}] $\rho$ induces a homeomorphism $Z/G_2\simeq
  G_1^{(0)}$.
\end{itemize}
\end{defi}

In Definition~\ref{def:generalized morphism}, one may replace (b)
by (b)' or (b)'' below:
\begin{itemize}
\item[{\rm{(b)'}}] $\rho$ is open and induces a bijection $Z/G_2\to
  G_1^{(0)}$.
\item[{\rm{(b)''}}] the map $Z\rtimes G_2\to Z\times_{G_1^{(0)}} Z$
  defined by $(z,\gamma)\mapsto (z,z\gamma)$ is
  a homeomorphism.
\end{itemize}

\begin{example}\label{ex:generalized morphism}
Let $G_1$ and $G_2$ be two groupoids.
If $f\colon G_1\to G_2$ is a groupoid morphism, let $Z=G_1^{(0)}
\times_{f,r}G_2$, $\rho(x,\gamma)=x$ and
$\sigma(x,\gamma)=s(\gamma)$.
Define the actions of $G_1$ and $G_2$ by $g\cdot(x,\gamma)\cdot\gamma'
=(r(g),f(g)\gamma\gamma')$. Then $(Z,\rho,\sigma)$ is a generalized
morphism from $G_1$ to $G_2$.
\end{example}

That $\rho$ is open follows from the fact that the range map $G_2\to
G_2^{(0)}$ is open and from Lemma~\ref{lem:open maps}.
The other properties in Definition~\ref{def:generalized morphism}
are easy to check.

\subsection{Locally proper generalized morphisms}
\begin{defi}
Let $G_1$ and $G_2$ be two groupoids
A generalized morphism
from $G_1$ to $G_2$ is said to be {\em locally proper} if the action of $G_1$
on $Z$ is $\sigma$-proper.
\end{defi}

Our terminology is justified by the following proposition:
\begin{prop}\label{prop:loc proper strict morphism}
Let $G_1$ and $G_2$ be two groupoids
such that $G_2^{(0)}$ is
Hausdorff. Let $f\colon G_1\to G_2$ be a groupoid morphism. Then
the associated generalized groupoid morphism is locally proper if and only
if the map $(f,r,s)\colon G_1\to G_2\times G_1^{(0)}\times G_1^{(0)}$
is proper.
\end{prop}

\pf
Let $\varphi\colon G_1\times_{f\circ s,r} G_2\to
(G_2\times_{s,s}G_2)\times_{r\times r,f\times f}
(G_1^{(0)}\times G_1^{(0)})$ defined by
$\varphi(g_1,g_2)=(f(g_1)g_2,g_2,r(g_1),s(g_1))$. By definition,
the action of
$G_1$ on $Z$ is proper if and only if $\varphi$ is a proper map.
Consider $\theta\colon G_2\times_{s,s}G_2 \to G_2^{(2)}$
given by
$(\gamma,\gamma')=(\gamma (\gamma')^{-1},\gamma')$. Let
$\psi=(\theta\times 1)\circ\varphi$. Since $\theta$ is a
homeomorphism, the action of $G_1$ on $Z$ is proper if and only if
$\psi$ is proper.

Suppose that $(f,r,s)$ is proper. Let $f'=(f,r,s)\times 1\colon
G_1\times G_2\to G_2\times G_1^{(0)}\times G_1^{(0)}\times G_2$.
Then $f'$ is proper.
Let $F=\{(\gamma,x,x',\gamma')\in G_2\times G_1^{(0)}\times
G_1^{(0)}\times G_2\vert\;
s(\gamma)=r(\gamma')=f(x'),\;r(\gamma)=f(x)\}$. Then $f'\colon
(f')^{-1}(F)\to F$ is proper, i.e. $\psi$ is proper.

Conversely, suppose that $\psi$ is proper. Let $F'=\{(\gamma,y,x,x')\in
G_2\times G_2^{(0)}\times G_1^{(0)} \times G_1^{(0)}\vert\;
s(\gamma)=y\}$. Then $\psi\colon \psi^{-1}(F')\to F'$ is proper, therefore
$(f,r,s)$ is proper.
\pfend

Our objective is now to show the
\begin{prop}\label{prop:composition of generalized morphisms}
Let $G_1$, $G_2$, $G_3$ be groupoids
Let
$(Z_1,\rho_1,\sigma_1)$ and $(Z_2,\rho_2,\sigma_2)$ be two generalized
groupoid morphisms from $G_1$ to $G_2$ and from $G_2$ to $G_3$
respectively. Then $(Z,\rho,\sigma)=(Z_1\times_{G_2} Z_2,\rho_1\times
1, 1\times \sigma_2)$ is a generalized groupoid morphism. If
$(Z_1,\rho_1,\sigma_1)$ and $(Z_2,\rho_2,\sigma_2)$ are locally proper, then
$(Z,\rho,\sigma)$ is locally proper. 
\end{prop}

Proposition~\ref{prop:composition of generalized morphisms} shows
that groupoids form a
category whose arrows are generalized morphisms, and that two
groupoids are isomorphic in that category if and only if they are
Morita-equivalent. Moreover, the same conclusions hold for the
category whose arrows are locally proper generalized morphisms. In
particular, local properness of generalized morphisms is invariant
under Morita-equivalence.

All the assertions of
Proposition~\ref{prop:composition of generalized morphisms}
follow from Lemma~\ref{lem:composition morphisms}.

\subsection{Proper generalized morphisms}
\begin{defi}\label{defi:proper gen morphism}
Let $G_1$ and $G_2$ be groupoids.
A generalized morphism $(Z,\rho,\sigma)$ from $G_1$ to $G_2$
is said to be {\em proper} if it is locally proper, and if for every
quasi-compact subspace $K$ of $G_2^{(0)}$, $\sigma^{-1}(K)$ is
$G_1$-compact.
\end{defi}

\begin{examples}\label{ex:proper-morphism}
\begin{itemize}
\item[{(a)}]
Let $X$ and $Y$ be locally compact spaces and $f\colon X\to Y$ a
continuous map. Then the generalized morphism
$(X,{\mathrm{Id}},f)$ is proper if and only if $f$ is proper.
\item[{(b)}]
Let $f\colon G_1\to G_2$ be a continuous morphism between
two locally compact groups. Let $p\colon G_2\to \{\ast\}$.
Then $(G_2,p,p)$ is proper if and only if $f$ is proper and
$f(G_1)$ is co-compact in $G_2$.
\item[{(c)}]
Let $G$ be a locally compact proper groupoid with Haar system such
that $G^{(0)}$ is Hausdorff, and let $\pi\colon G^{(0)}\to G^{(0)}/G$
be the canonical mapping. Then $(G^{(0)},{\mathrm{Id}},\pi)$ is
a proper generalized morphism from $G$ to $G^{(0)}/G$.
\end{itemize}
\end{examples}

\subsection{Construction of a $C^*$-correspondence}
Until the end of the section, our goal is to prove:
\begin{theo}\label{thm:correspondence}
Let $G_1$ and $G_2$ be locally compact
groupoids with Haar system
such that $G_1^{(0)}$ and $G_2^{(0)}$ are Hausdorff, and
$(Z,\rho,\sigma)$ a locally
proper generalized morphism from $G_1$ to $G_2$.
Then one can construct a $C^*_r(G_1)$-Hilbert module
${\mathcal{E}}_Z$ and a map $\pi\colon C^*_r(G_2)\to
{\mathcal{L}}({\mathcal{E}}_Z)$.
Moreover, if $(Z,\rho,\sigma)$ is proper, then
$\pi$ maps to ${\mathcal{K}}({\mathcal{E}}_Z)$. Therefore, it
gives an element of $KK(C^*_r(G_2),C^*_r(G_1))$.
\end{theo}


\begin{coro}(see \cite{mrw})
Let $G_1$ and $G_2$ be locally compact
groupoids with Haar system
such that $G_1^{(0)}$ and $G_2^{(0)}$ are Hausdorff. If
$G_1$ and $G_2$ are Morita-equivalent, then $C^*_r(G_1)$ and
$C^*_r(G_2)$ are Morita-equivalent.
\end{coro}

\begin{coro}
Let $f\colon G_1\to G_2$ be morphism between two locally compact
groupoids with Haar system
such that $G_1^{(0)}$ and $G_2^{(0)}$ are Hausdorff.
If the restriction of $f$ to $(G_1)_K^K$ is proper
for each compact set $K\subset (G_1)^{(0)}$ then $f$ induces
a correspondence ${\mathcal{E}}_f$ from $C^*_r(G_2)$ to
$C^*_r(G_1)$. If in addition for every compact set $K\subset G_2^{(0)}$
the quotient of $G_1^{(0)}\times_{f,r}(G_2)_K$ by the diagonal
action of $G_1$ is compact, then $C^*_r(G_2)$ maps to
${\mathcal{K}}({\mathcal{E}}_f)$ and thus $f$ defines a $KK$-element
$[f]\in KK(C^*_r(G_2),C^*_r(G_1))$.
\end{coro}

\pf
See Proposition~\ref{prop:loc proper strict morphism} and
Definition~\ref{defi:proper gen morphism} applied
to the generalized morphism $Z_f=G_1^{(0)}\times_{f,r}G_2$
as in Example~\ref{ex:generalized morphism}
\pfend


The rest of the section is devoted to proving
Theorem~\ref{thm:correspondence}.

Let us first recall the construction of the correspondence when
the groupoids are Hausdorff \cite{mo}. It is the closure of
$C_c(Z)$ with the $C^*_r(G_1)$-valued scalar product
\begin{equation}\label{eqn:scalar-product}
\langle\xi,\eta\rangle(g) = \int_{\gamma\in (G_2)^{\sigma(z)}}
\overline{\xi(z\gamma)}
\eta(g^{-1}z\gamma)\,\lambda^{\sigma(z)}({\mathrm{d}}\gamma),
\end{equation}
where $z$ is an arbitrary element of $Z$ such that $\rho(z)=r(g)$.
The right $C^*_r(G_1)$-module structure is defined
$\forall \xi\in C_c(Z)$, $\forall a\in C_c(G_1)$ by
\begin{equation}\label{eqn:right-module}
(\xi a)(z)=\int_{g\in (G_1)^{\rho(z)}}
\xi(g^{-1}z)a(g^{-1})\,\lambda^{\rho(z)}({\mathrm{d}}g),
\end{equation}
and the left action of $C^*_r(G_2)$ is
\begin{equation}\label{eqn:left-action}
(b\xi)(z)=\int_{\gamma\in (G_2)^{\sigma(z)}}
  b(\gamma)\xi(z\gamma)\,\lambda^{\sigma(z)}({\mathrm{d}}\gamma)
\end{equation}
for all $b\in C_c(G_2)$.
\par\medskip

We now come back to non-Hausdorff groupoids.
For every open Hausdorff set $V\subset Z$, denote by $V'$ its closure
in ${\mathcal{H}}((G_1\ltimes Z)_V^V)$,
where $z\in V$ is identified to $(\rho(z),z)\in
{\mathcal{H}}((G_1\ltimes Z)_V^V)$.
Let ${\mathcal{E}}^0_V$ be the
set of $\xi\in C_c(V')$ such that
$\displaystyle\xi(z)=\frac{\xi(S\times\{z\})}{\sqrt{\# S}}$ for all
$S\times\{z\}\in V'$.

\begin{lem}\label{lem:E_Z independent of cover}
The space ${\mathcal{E}}_Z^0=\sum_{i\in I}{\mathcal{E}}^0_{V_i}$ is
independent of the choice of the cover $(V_i)$ of $Z$ by Hausdorff
open subspaces.
\end{lem}

\pf
It suffices to show that for every open Hausdorff subspace $V$ of $Z$,
one has ${\mathcal{E}}^0_V\subset
\sum_{i\in I}{\mathcal{E}}^0_{V_i}$. Let
$\xi\in {\mathcal{E}}^0_V$. Denote by $q_V\colon V'\to V$ the
canonical map defined by $q_V(S\times\{z\})=z$. Let $K\subset V$
compact such that ${\mathrm{supp}}(\xi)\subset q_V^{-1}(K)$. There
exists $J\subset I$ finite such that $K\subset \cup_{j\in J}
V_j$. Let $(\varphi_j)_{j\in J}$ be a partition of unity associated
to that cover, and $\xi_j=\xi . (\varphi_j\circ q_V)$. One easily
checks that $\xi_j\in {\mathcal{E}}^0_{V_j}$ and that $\xi=\sum_{j\in
  J}\xi_j$.
\pfend

We now define a $C^*_r(G_1)$-valued scalar product on
${\mathcal{E}}^0_Z$ by Eqn. (\ref{eqn:scalar-product})
where $z$ is an arbitrary element of $Z$ such that $\rho(z)=r(g)$. Our
definition is independent of the choice of $z$, since if $z'$ is
another element, there exists $\gamma'\in G_2$ such that
$z'=z\gamma'$, and the Haar system on $G_2$ is left-invariant.

Moreover, the integral is convergent for all $g\in G_1$ because the
action of $G_2$ on $Z$ is proper.

Let us show that $\langle\xi,\eta\rangle\in C_c(G_1)$ for all $\xi$,
$\eta\in {\mathcal{E}}^0_Z$. We need a preliminary lemma:

\begin{lem}\label{lem:limit proper map}
Let $X$ and $Y$ be two topological spaces such that $X$ is locally
compact and $f\colon X\to Y$
proper. Let ${\mathcal{F}}$ be a ultrafilter such that
$f$ converges to $y\in Y$ with respect to ${\mathcal{F}}$. Then there exists
$x\in X$ such that $f(x)=y$ and ${\mathcal{F}}$ converges to $x$.
\end{lem}

\pf
Let $Q=f^{-1}(y)$. Since $f$ is proper, $Q$ is quasi-compact. Suppose
that for all $x\in Q$, ${\mathcal{F}}$ does not converge to $x$.
Then there exists an open neighborhood $V_x$ of $x$ such that
$V_x^c\in {\mathcal{F}}$. Extracting a finite cover $(V_1,\ldots,V_n)$
of $Q$, there exists an open neighborhood $V$ of $Q$ such that
$V^c\in {\mathcal{F}}$. Since $f$ is closed, $f(V^c)^c$ is a
neighborhood of $y$. By assumption, $f(V^c)^c\in f({\mathcal{F}})$,
i.e. $\exists A\in {\mathcal{F}}$, $f(A)\subset f(V^c)^c$.
This implies that $A\subset V$, therefore $V\in{\mathcal{F}}$:
this contradicts $V^c\in {\mathcal{F}}$.

Consequently, there exists $x\in Q$ such that ${\mathcal{F}}$
converges to $x$.
\pfend

To show that $\langle\xi,\eta\rangle\in C_c(G_1)$,
we can suppose that $\xi\in
{\mathcal{E}}^0_U$ and $\eta\in {\mathcal{E}}^0_V$, where $U$ and $V$
are open Hausdorff. Let $F(g,z)=\overline{\xi(z)}\eta(g^{-1}z)$,
defined on $\Gamma=G_1\times_{r,\rho} Z$.
Since the action of $G_1$ on $Z$ is proper, $F$ is quasi-compactly
supported. Let us show that $F\in C_c(\Gamma)$.

Let ${\mathcal{F}}$ be a ultrafilter on $\Gamma$, convergent in
${\mathcal{H}}\Gamma$. Since $G_1^{(0)}$ is Hausdorff,
its limit has the form $S=S'g_0\times S''$
where $S'\subset (G_1)_{r(g_0)}^{r(g_0)}$, $S''\subset \rho^{-1}(r(g_0))$.
Moreover, $S'$ is a subgroup of $(G_1)_{r(g)}^{r(g)}$ by the
proof of Lemma~\ref{lem:S is a group}.

Suppose that there exist $z_0,z_1\in S''$ and $g_1\in S'g_0$ such that
$z_0\in U$ and $g_1^{-1}z_1\in V$. By Lemma~\ref{lem:limit proper map}
applied to the proper map $G_1\rtimes Z\to Z\times Z$, there exists
$s_0\in S'$ such that $z_0=s_0z_1$. We may assume that $g_0=s_0g_1$.
Then $\sum_{s\in S} F(s)=\sum_{s'\in S'}\overline{\xi(z_0)}
\eta(g_0^{-1}(s')^{-1} z_0)$. If $s'\notin{\mathrm{stab}}(z_0)$, then
$g_0^{-1}(s')^{-1}z_0\notin V$ since $g_0^{-1}z_0$ and
$g_0^{-1}(s')^{-1} z_0$ are distinct limits of
$(g,z)\mapsto g^{-1}z$ with respect to
${\mathcal{F}}$ and $V$ is Hausdorff. Therefore,
\begin{eqnarray*}
\sum_{s\in S} F(s)&=&\#({\mathrm{stab}}(z_0)\cap S')
\overline{\xi(z_0)} \eta(g_0^{-1}z_0)\\
&=&\overline{\sqrt{\#({\mathrm{stab}}(z_0)\cap S')}\xi(z_0)}
\sqrt{\#({\mathrm{stab}}(g_0^{-1}z_0)\cap (g_0^{-1}S'g_0))}\eta(z_0)\\
&=&\lim_{{\mathcal{F}}} \overline{\xi(z)}\eta(g^{-1}z) =
\lim_{{\mathcal{F}}} F(g,z).
\end{eqnarray*}

If for all $z_0$, $z_1\in S''$ and all $g_1\in S'g_0$,
$(z_0,g_1^{-1}z_1)\notin U\times V$, then
$\sum_{s\in S} F(g,z)=0=\lim_{{\mathcal{F}}} F(g,z)$.

By Proposition~\ref{prop:C_c(X)}, $F\in C_c(\Gamma)$.

Since $\langle \xi,\eta\rangle (g) = \int_{\gamma\in
  (G_2)^{\sigma(z)}} F(g,z\gamma)\,\lambda^{\sigma(z)}({\mathrm{d}}\gamma)$, to
  prove that $\langle\xi,\eta\rangle\in C_c(G_1)$ it
  suffices to show:

\begin{lem}\label{lem:belongs to Cc}
Let $G_1$ and $G_2$ be two locally compact groupoids with Haar system
such that $G_i^{(0)}$ are Hausdorff.
Let $(Z,\rho,\sigma)$ be a generalized morphism from $G_1$
to $G_2$. Let $\Gamma=G_1\times_{r,\rho} Z$.
Then for every $F\in C_c(\Gamma)$, the function
$$g\mapsto \int_{\gamma\in (G_2)^{\sigma(z)}}
F(g,z\gamma)\,\lambda^{\sigma(z)}({\mathrm{d}}\gamma),$$
where $z\in Z$ is an arbitrary element such that $\rho(z)=r(g)$,
belongs to $C_c(G_1)$.
\end{lem}
 
\pf
Suppose first that $F(g,z)= f(g)h(z)$, where $f\in C_c(G_1)$ and $h\in
C_c(Z)$. Let $H(z)=\int_{\gamma\in (G_2)^{\sigma(z)}} h(z\gamma)\,
\lambda^{\sigma(z)}({\mathrm{d}}\gamma)$. By Lemma~\ref{lem:f(s(g))} below
(applied to the groupoid $Z\rtimes G_2$), $H$ is
continuous. It is obviously $G_2$-invariant, therefore $H\in
C_c(Z/G_2)$. Let $\tilde H\in C_c(G_1^{(0)})\simeq C_c(Z/G_2)$
correspond to $H$.
The map
$$g\mapsto \int_{\gamma\in (G_2)^{\sigma(z)}}
F(g,z\gamma)\,\lambda^{\sigma(z)}({\mathrm{d}}\gamma) = f(g){\tilde
  H}(s(g))$$
thus belongs to $C_c(G_1)$.

By linearity, the lemma is true for $F\in C_c(G_1)\otimes C_c(Z)$.
By Lemma~\ref{lem:density C_c(X)xC_c(Y)} and
Lemma~\ref{lem:C_c(X) to C_c(Y)},  $F$ is the uniform limit of
functions $F_n\in C_c(G_1)\otimes C_c(Z)$ which are supported in a
fixed quasi-compact set $Q=Q_1\times Q_2\subset G_1\times Z$. Let
$Q'\subset Z$ quasi-compact such that $\rho(Q')\supset r(Q_1)$. Since
the action of $G_2$ on $Z$ is proper, $K=\{\gamma\in G_2\vert\;
Q'\gamma\cap Q_2\ne\emptyset\}$ is quasi-compact. Using
the fact that $G_1^{(0)}\simeq Z/G_2$, it is easy to see that
$${\sup_{(g,z)\in\Gamma}\int_{\gamma\in (G_2)^{\sigma(z)}}
1_Q(g,z\gamma)\,\lambda^{\sigma(z)}({\mathrm{d}}\gamma)
\le\sup_{z\in Q'} \int_{\gamma\in G_2^{\sigma(z)}}
1_{Q_2}(z\gamma)\lambda^{\sigma(z)}({\mathrm{d}}\gamma)}$$
$$\le \sup_{x\in G_2^{(0)}} \int_{\gamma\in G_2^x} 1_K(\gamma)
\lambda^{x}({\mathrm{d}}\gamma)<\infty$$

by Lemma~\ref{lem:measure of quasi-compact finite}. Therefore,
$$\lim_{n\to\infty}\sup_{g\in G_1}
\left|
\int_{\gamma\in G_2^{\sigma(z)}} F(g,z\gamma)-F_n(g,z\gamma)\,
\lambda^{\sigma(z)}({\mathrm{d}}\gamma)
\right|
=0.$$
The conclusion follows from Corollary~\ref{cor:limit in C_c}.
\pfend

In the proof of Lemma~\ref{lem:belongs to Cc} we used the

\begin{lem}\label{lem:f(s(g))}
Let $G$ be a locally compact, proper groupoid with Haar system, such
that $G^x$ is Hausdorff for all $x\in G^{(0)}$, and $G_x^x=\{x\}$ for
all $x\in G^{(0)}$. We do not assume $G^{(0)}$ to be
Hausdorff. Then $\forall f\in C_c(G^{(0)})$,
$$\varphi\colon G^{(0)}\to\C,\quad
x\mapsto \int_{g\in G^x} f(s(g))\,\lambda^x(dg)$$
is continuous.
\end{lem}

\pf
Let $V$ be an open, Hausdorff subspace of $G^{(0)}$. Let $h\in
C_c(V)$. Since $(r,s)\colon
G\to G^{(0)}\times G^{(0)}$ is
a homeomorphism from $G$ onto a closed subspace of $G^{(0)}\times
G^{(0)}$, and $(x,y)\mapsto h(x)f(y)$ belongs to $C_c(G^{(0)}\times
G^{(0)})$, the map $g\mapsto h(r(g))f(s(g))$ belongs to $C_c(G)$,
therefore by definition of a Haar system, $x\mapsto \int_{g\in G^x}
h(r(g))f(s(g))\,\lambda^x(dg) = h(x)\varphi(x)$ belongs to
$C_c(G^{(0)})$.

Since $h\in C_c(V)$ is arbitrary, this shows that $\varphi_{\vert V}$
is continuous, hence $\varphi$ is continuous on $G^{(0)}$.
\pfend

Now, let us show the positivity of the scalar product.
Recall that for all $x\in G_1^{(0)}$ there is a representation
$\pi_{G_1,x}\colon C^*(G_1)\to {\mathcal{L}}(L^2(G_1^x))$ such that
for all $a\in C_c(G_1)$ and all $\eta\in C_c(G_1^x)$,
$$(\pi_{G_1,x}(a) \eta)(g) = \int_{h\in G_1^{s(g)}}
a(h)\eta(gh)\,\lambda^{s(g)}({\mathrm{d}}h).$$
By definition, $\|a\|_{C^*_r(G_1)}=\sup_{x\in G_1^{(0)}}
\|\pi_{G_1,x}(a)\|$.

\begin{eqnarray*}
\langle\eta,\pi_{G_1,x}(a)\eta\rangle &=&
\int_{g\in G_1^x,\;h\in G_1^{s(g)}}\overline{\eta(g)}
a(h)\eta(gh)\,\lambda^{s(g)}({\mathrm{d}}h)\lambda^x(dg)\\
&=&\int_{g\in G_1^x,\;h\in G^{s(g)}}\overline{\eta(g)}a(g^{-1}h)\eta(h)
\,\lambda^x(dg)\lambda^x(dh).
\end{eqnarray*}
Fix $z\in Z$ such that $\rho(z)=x$. Replacing $a(g^{-1}h)$ by
$$\langle\xi,\xi\rangle (g^{-1}h)
=\int_{\gamma\in G_2^{\sigma(z)}}
\overline{\xi(g^{-1}z\gamma)}\xi(h^{-1}z\gamma)\,
\lambda^{\sigma(z)}({\mathrm{d}}\gamma),$$
we get
\begin{equation}\label{eqn:scalar product}
{\langle\eta,\pi_{G_1,x}(\langle\xi,\xi\rangle)\eta\rangle
= \int_{\gamma\in G_2^{\sigma(z)}} \lambda^{\sigma(z)}({\mathrm{d}}\gamma)}
\left|
\int_{g\in G^x}\eta(g)\xi(g^{-1}z\gamma)\,\lambda^x(dg)
\right|^2.
\end{equation}
It follows that $\pi_{G_1,x}(\langle\xi,\xi\rangle)\ge 0$ for all $x\in
G_1^{(0)}$, so $\langle\xi,\xi\rangle\ge 0$ in $C^*_r(G_1)$.
\par\medskip

Now, let us define a
$C^*_r(G_1)$-module structure on ${\mathcal{E}}_Z^0$ by
Eqn.(\ref{eqn:right-module})
%
for all $\xi\in {\mathcal{E}}^0_Z$ and $a\in C_c(G_1)$.

Let us show that $\xi a\in {\mathcal{E}}^0_Z$.
We need a preliminary lemma:

\begin{lem}\label{lem:product compact}
Let $X$ and $Y$ be quasi-compact spaces, $(\Omega_k)$ an open cover of
$X\times Y$. Then there exist finite open covers $(X_i)$ and $(Y_j)$
of $X$ and $Y$ such that $\forall i,j$ $\exists k$, $X_i\times
Y_j\subset \Omega_k$.
\end{lem}

\pf
For all $(x,y)\in X\times Y$ choose open neighborhoods $U_{x,y}$ and
$V_{x,y}$ of $x$ and $y$ such that $U_{x,y}\times V_{x,y}\subset
\Omega_k$ for some $k$. For $y$ fixed, there exist $x_1,\ldots,x_n$
such that $(U_{x_i,y})_{1\le i\le n}$ covers $X$. Let
$V_y=\cap_{i=1}^n U_{x_i,y}$. Then for all $(x,y)\in X\times Y$, there
exists an open neighborhood $U'_{x,y}$ of $x$ and $k$ such that
$U'_{x,y}\times V_y\subset \Omega_k$.

Let $(V_1,\ldots,V_m)=(V_{y_1},\ldots,V_{y_m})$ such that $\cup_{1\le
  j\le m} V_j=Y$. For all $x\in X$, let $U'_x=\cap_{j=1}^m
  U'_{x,y_j}$. Let $(U_1,\ldots,U_p)$ be a finite sub-cover of
  $(U'_x)_{x\in X}$. Then for all $i$ and for all $j$, there exists
  $k$ such that $U_i\times V_j\subset\Omega_k$.
\pfend

Let $Q_1$ and $Q_2$ be quasi-compact subspaces of $G_1$ of $Z$
respectively such that $a^{-1}(\C^*)\subset {Q}_1$ and
$\xi^{-1}(\C^*)\subset {Q}_2$.
Let $Q$ be a quasi-compact subspace of $Z$ such that $\forall g\in
Q_1$, $\forall z\in Q_2$, $g^{-1}z\in Q$. Let $(U_k)$ be a finite
cover of $Q$ by Hausdorff open subspaces of $Z$. Let
$Q'=Q_1\times_{r,\rho} Q_2$. Then $Q'$ is a closed subspace of
$Q_1\times Q_2$. Let
$\Omega'_k=\{(g,z)\in Q'\vert\; g^{-1}z\in
U_k\}$. Then $(\Omega'_k)$ is a finite open cover of
$Q'$. Let $\Omega_k$ be an open subspace of
$Q_1\times Q_2$ such that $\Omega'_k=\Omega_k\cap
Q'$. Then $\{Q_1\times Q_2-Q'\} \cup \{\Omega_k\}$ is an open
cover of $Q_1\times Q_2$. Using Lemma~\ref{lem:product compact}, there
exist finite families of Hausdorff open sets $(W_i)$ and $(V_j)$
which cover $Q_1$ and $Q_2$, such that for
all $i$, $j$ and for all $(g,z)\in W_i\times_{G_1^{(0)}} V_j$, there
exists $k$ such that $g^{-1}z\in U_k$.

Thus, we can assume by linearity and by
Lemmas~\ref{lem:predetermined cover}
and~\ref{lem:E_Z independent of cover} that $\xi\in {\mathcal{E}}^0_V$,
$a\in C_c(W)$, $U=W^{-1}V$, and $U$, $V$ and $W$ are open and
Hausdorff.

Let $\Omega=\{(g,S)\in W^{-1}\times U'\vert\; g^{-1} q_U(S)\in V\}$. Then
the map $(g,S)\mapsto (g^{-1},g^{-1}S)$ is a homeomorphism from
$\Omega$ onto $W\times_{r,\rho\circ q_V}V'$. Therefore, the map
$(g,z)\mapsto \xi(g^{-1}z)a(g^{-1})$ belongs to $C_c(\Omega)\subset
C_c(G_1\times_{r,\rho\circ q_V}U')$. By Lemma~\ref{lem:Haar product},
$$S\mapsto (\xi a)(S)=\int_{g\in G_1^{\rho\circ q_V(S)}}
  \xi(g^{-1}S)a(g^{-1})   \,\lambda^{\rho\circ q_V(S)}({\mathrm{d}}g)$$
belongs
  to $C_c(U')$. It is immediate that $(\xi a)(S)=\sqrt{\# S}(\xi
  a)( q(S))$ for all $S\in U'$, therefore $\xi a\in
  {\mathcal{E}}^0_U$.
This completes the proof that $\xi a\in {\mathcal{E}}^0_Z$.
\par\medskip

Finally, it is not hard to check that $\langle\xi,\eta a\rangle
=\langle\xi,\eta\rangle\ast a$. Therefore, the completion
${\mathcal{E}}_Z$ of ${\mathcal{E}}^0_Z$ with respect to the norm
$\|\xi\|=\|\langle\xi,\xi\rangle\|^{1/2}$ is a $C^*_r(G_1)$-Hilbert
module.
\par\medskip

Let us now construct a morphism $\pi\colon C^*_r(G_2)\to
{\mathcal{L}}({\mathcal{E}}_Z)$. For every $\xi\in {\mathcal{E}}_Z^0$
and every $b\in C_c(G_2)$, define $b\xi$ by Eqn.(\ref{eqn:left-action}).
Let us check that $b\xi\in {\mathcal{E}}^0_Z$. As above,
by linearity we may
assume that $\xi\in {\mathcal{E}}^0_V$, $b\in C_c(W)$ and
$VW^{-1}\subset U$, where $V\subset Z$, $U\subset Z$ and $W\subset
G_2$ are open and Hausdorff.

Let $\Phi(S,\gamma)=(S\gamma,\gamma)$. Then $\Phi$ is a homeomorphism
from $\Omega=\{(S,\gamma)\in U'\times_{\sigma\circ q_U,r}W\vert\;
 q_U(S)\gamma\in V\}$ onto $V'\times_{\sigma\circ q_V,s}W$. Let
$F(z,\gamma)=b(\gamma)\xi(z\gamma)$. Since $F=(\xi\otimes b)\circ
\Phi$, $F$ is an element of $C_c(\Omega)\subset
C_c(U'\times_{\sigma\circ q_U,r}W)$.
By Lemma~\ref{lem:Haar product}, $b\xi\in C_c(U')$.

It is immediate that $(b\xi)(S)=\sqrt{\# S}
(b\xi)( q(S))$. Therefore, $b\xi\in {\mathcal{E}}^0_U\subset
{\mathcal{E}}^0_Z$.

Let us prove that $\|b\xi\|\le \|b\| \, \|\xi\|$.
Let
$$\zeta(\gamma)=\int_{g\in G_1^x} \eta(g)\xi(g^{-1}z\gamma)\,
\lambda^x(dg),$$
 where $z\in Z$ such that $\rho(z)=r(g)$ is arbitrary.
From~(\ref{eqn:scalar product}),
$$\langle \eta,\pi_{G_1,x}(\langle\xi,\xi\rangle)\eta\rangle
=\|\zeta\|^2_{L^2(G_2^{\sigma(z)})}.$$
A similar calculation shows that
$$\langle\eta,\pi_{G_1,x}(\langle b\xi,b\xi\rangle)\eta\rangle
=\int_{\gamma\in G_2^{\sigma(z)}}\!\!\!\!\!\!\!\!
\lambda^{\sigma(z)}({\mathrm{d}}\gamma)
\left|\int_{g\in G_1^x}\!\!\!\!\!\eta(g)\xi(g^{-1}z\gamma\gamma')
b(\gamma') \,
  \lambda^{s(\gamma)}({\mathrm{d}}\gamma')\right|^2$$
$$=\langle b\zeta,b\zeta\rangle \le \|b\|^2\|\zeta\|^2.$$

By density of $C_c(G_2^x)$ in $L^2(G_2^x)$, $\|\pi_{G_1,x}(\langle
b\xi, b\xi\rangle)\| \le \|b\|^2
\|\pi_{G_1,x}(\langle\xi,\xi\rangle)\|$.
Taking the supremum over $x\in G_1^{(0)}$, we get
$\|b\xi\|\le \|b\|\,
\|\xi\|$. It follows that $b\mapsto (\xi\mapsto b\xi)$ extends to a
$*$-morphism $\pi\colon
C^*_r(G_2)\to {\mathcal{L}}({\mathcal{E}}_Z)$.
\par\medskip

Finally, suppose now that $(Z,\rho,\sigma)$ is proper, and
let us show that $C^*_r(G_2)$ maps to ${\mathcal{K}}
({\mathcal{E}}_Z)$.

For every $\eta$, $\zeta\in {\mathcal{E}}_Z^0$, denote by
$T_{\eta,\zeta}$ the operator $T_{\eta,\zeta}(\xi) = \eta\langle
\zeta,\xi\rangle$. Compact operators are elements of the closed linear
span of $T_{\eta,\zeta}$'s. Let us write an explicit formula for
$T_{\eta,\zeta}$:
\begin{eqnarray*}
T_{\eta,\zeta}(\xi)(z) &=&
\int_{g\in G_1^{\rho(z)}} \eta(g^{-1}z)\langle \zeta,\xi\rangle
(g^{-1})\,\lambda^{\rho(z)}({\mathrm{d}}g)\\
&=&\int_{g\in G_1^{\rho(z)}} \eta(g^{-1}z) \int_{\gamma\in
  G_2^{\sigma(z)}} \overline{\zeta(g^{-1}z\gamma)} \xi(z\gamma)\,
\lambda^{\sigma(z)}({\mathrm{d}}\gamma)\lambda^{\rho(z)}({\mathrm{d}}g).
\end{eqnarray*}

Let $b\in C_c(G_2)$, let us show that
$\pi(b)\in{\mathcal{K}}({\mathcal{E}}_Z)$. Let $K$ be a quasi-compact
subspace of $G_2$ such that $b^{-1}(\C^*)\subset K$.
Since $(Z,\rho,\sigma)$
is a proper generalized morphism, there exists a quasi-compact
subspace $Q$ of $Z$ such that $\sigma^{-1}(r(K))\subset G_1\rond{Q}$.
Before we proceed, we need a lemma:

\begin{lem}\label{lem:Omega}
Let $G_2$ be a locally compact groupoid acting freely and properly on
a locally compact space $Z$ with momentum map $\sigma\colon Z\to G_2^{(0)}$.
Then for every $(z_0,\gamma_0)\in
Z\rtimes G_2$, there exists a Hausdorff open neighborhood
$\Omega_{z_0,\gamma_0}$ of $(z_0,\gamma_0)$ such that
\begin{itemize}
\item $U=\{z_1\gamma_1\vert\;(z_1,\gamma_1)
\in \Omega_{z_0,\gamma_0}\}$ is Hausdorff;
\item there exists a Hausdorff open neighborhood $W$ of $\gamma_0$
such that $\forall \gamma\in G_2$,
$\forall z\in pr_1(\Omega_{z_0,\gamma_0})$,
$\forall z'\in U$, $z'=z\gamma \implies \gamma\in W$.
\end{itemize}
\end{lem}

\pf
Let $R=\{(z,z')\in Z\times Z\vert\; \exists\gamma\in G_2,\; z'=z\gamma\}$.
Since the $G_2$-action is free and proper, there exists a
continuous function $\phi\colon R\to G_2$ such that
$\phi(z,z\gamma)=\gamma$.
Let $W$ be an open Hausdorff neighborhood
of $\gamma_0$. By continuity of $\phi$, there exist
open Hausdorff neighborhoods $V$ and $U_0$ of $z_0$ and
$z_0\gamma_0$ such that for all $(z,z')\in
R\cap(V\times U_0)$, $\phi(z,z')\in W$. By continuity of the action,
there exists an open neighborhood $\Omega_{z_0,\gamma_0}$ of
$(z_0,\gamma_0)$  such that $\forall (z_1,\gamma_1)\in
\Omega_{z_0,\gamma_0}$, $z_1\gamma_1\in U_0$ and $z_1\in V$.
\pfend

By Lemma~\ref{lem:product compact}, there exist finite covers $(V_i)$
of $Q$ and $(W_j)$ of $K$ such that for every $i$, $j$,
$(Z\times_{G_2^{(0)}}G_2)\cap (V_i\times W_j)\subset
\Omega_{z_0,\gamma_0}$ for some $(z_0,\gamma_0)$.

By Lemma~\ref{lem:local construction of c} applied to the groupoid
$(G_1\ltimes Z)_{V_i}^{V_i}$, for all $i$ there exists
$c'_i\in C_c(V'_i)_+$ such that
$c'_i(S)={(\#S)} c'_i( q_{V_i}(S))$ for all $S\in V'_i$, and such that
$\sum_i c'_i \ge 1$ on $Q$. Let
$$f_i(z)=\int_{g\in G_1^{\rho(z)}}
c'_i(g^{-1}z) \,\lambda^{\rho(z)}({\mathrm{d}}g)$$
and let $f=\sum_i f_i$.
As in the proof of Theorem~\ref{thm:cutoff}, one can show that
for every Hausdorff open subspace $V$ of $Z$ and every $h\in C_c(V)$,
$(g,z)\mapsto h(z)c'_i(g^{-1}z)$ belongs to $C_c(G\ltimes Z)$, therefore
$hf_i$ is continuous on $V$. Since $h$ is arbitrary, it follows that
$f_i$ is continuous, thus $f$ is continuous.
Moreover, $f$ is
$G_1$-equivariant, nonnegative, and $\inf_{Q}f>0$. Therefore, there
exists $f_1\in C_c(G_1\backslash Z)$ such that $f_1(z)=1/f(z)$ for all
$z\in Q$. Let $c_i(z)=f_1(z)c'_i(z)$. Let
$$T_i(\xi)(z) = \int_{g\in
  G_1^{\rho(z)}} \int_{\gamma\in G_2^{\sigma(z)}}
c_i(g^{-1}z)b(\gamma)\xi(z\gamma) \,\lambda^{\rho(z)}({\mathrm{d}}g)
\lambda^{\sigma(z)}({\mathrm{d}}\gamma).$$
Then $\pi(b)=\sum_i T_i$, therefore it
suffices to show that $T_i$ is a compact operator for all $i$.

By linearity and by Lemma~\ref{lem:predetermined cover},
one may assume that $b\in C_c(W_j)$ for some $j$. Then,
by construction of $V_i$ (see Lemma~\ref{lem:Omega}),
there exist open Hausdorff sets $U\subset Z$ and $W\subset G_2$ such
that $\{\gamma\in G_2\vert\;\exists (z,z')\in V_i\times U,\;
z'=z\gamma\} \subset W$, and $\{z\gamma\vert\;
(z,\gamma)\in V_i\times_{\sigma,r} W\} \subset U$.

The map $(z,z\gamma)\mapsto c(z)b(\gamma)$ defines an element of
$C_c(V'_i\times U)$. Let $L_1\times L_2\subset V_i\times U$ compact
such that $(z,z\gamma)\mapsto c(z)b(\gamma)$ is supported on
$ q_{V_i}^{-1}(L_1) \times L_2$.
By Lemma~\ref{lem:local construction of c} applied to the groupoids
$(G_1\ltimes Z)_{V_i}^{V_i}$ and $(G_1\ltimes Z)_U^U$, there exist
$d_1\in C_c(V'_i)_+$ and $d_2\in
C_c(U')_+$ such that $d_1>0$ on $L_1$ and $d_2>0$ on $L_2$,
$d_1(S)=\sqrt{\# S} d_1( q_{V_i}(S))$ for all $S\in V'_i$, and
$d_2(S)=\sqrt{\# S} d_2( q_U(S))$ for all $S\in U'$.
Let
$$f(z,z\gamma)=\frac{c(z)b(\gamma)}{d_1(z)d_2(z\gamma)}.$$
Then $f\in C_c(V_i\times_{G_1^{(0)}} U)$. Therefore, $f$ is the
uniform limit of a sequence $f_n=\sum \alpha_{n,k}\otimes
\overline{\beta_{n,k}}$ in $C_c(V_i)\otimes C_c(U)$ such that all the
$f_n$ are supported in a fixed compact set. Then $T_i$ is the
norm-limit of
$\sum_k T_{d_1\alpha_{n,k},d_2\beta_{n,k}}$, therefore it is compact.

\begin{rem}
The construction in Theorem~\ref{thm:correspondence} is functorial
with respect to the composition of generalized morphisms and
of correspondences. We don't include a proof of this fact, as
it is tedious but elementary. It is an easy exercise when
$G_1$ and $G_2$ are Hausdorff.
\end{rem}

\bibliographystyle{plain}

\Addresses

\end{document}